\documentclass[11pt,reqno]{amsart}
\usepackage[T1]{fontenc}
\usepackage{lmodern}
\usepackage{amsmath,amssymb,mathtools,mathrsfs}
\usepackage{microtype}
\usepackage{xcolor}
\usepackage{enumitem}
\usepackage[a4paper,margin=28mm]{geometry}
\usepackage[colorlinks=true,linkcolor=blue,citecolor=blue,urlcolor=blue]{hyperref}
\numberwithin{equation}{section}
\newtheorem{theorem}{Theorem}[section]
\newtheorem{lemma}[theorem]{Lemma}
\newtheorem{hypothesis}[theorem]{Hypothesis}
\newtheorem{proposition}[theorem]{Proposition}
\newtheorem{corollary}[theorem]{Corollary}
\theoremstyle{definition}
\newtheorem{definition}[theorem]{Definition}
\newtheorem{example}[theorem]{Example}
\theoremstyle{remark}
\newtheorem{remark}[theorem]{Remark}
\newtheorem*{remark*}{Remark}

\newcommand{\dd}{\,\mathrm{d}}
\newcommand{\norm}[1]{\lVert #1\rVert}
\newcommand{\pair}[2]{\langle #1,#2\rangle}
\setlist[enumerate]{leftmargin=*,label=\textup{(\roman*)}}
\allowdisplaybreaks[1]
\title[Uniform stabilization of degenerate thermoelasticity]
{Uniform stabilization of degenerate thermoelasticity
with non-autonomous conductivities}
\author[M. Liao]{Menglan Liao}
\address{School of Mathematics, Hohai University, Nanjing, 210098, China \&
Laboratory of Mathematical Modeling and Intelligent Computing for Water Systems, Hohai University, Nanjing, 211100, China.}
\email{liaoml@hhu.edu.cn}
\thanks{This work was supported by the National Natural Science Foundation of China (12401290) and the Natural Science Foundation of Jiangsu Province (BK20230946).}
\subjclass[2020]{Primary 35B40; Secondary 35K65, 35L05, 47D06.}
\keywords{Degenerate thermoelasticity, non-autonomous evolution equations, global well-posedness, energy decay, multiplier method, Lyapunov functional}
\date{}

\begin{document}
\begin{abstract}
We study a one-dimensional thermoelastic system with space- and time-dependent conductivity vanishing at one endpoint, under Dirichlet or zero-flux thermal conditions. Uniform comparison with the initial conductivity and local temporal regularity yield global well-posedness for energy data. A quadratic lower bound ensures uniform exponential decay for every fixed coupling $\kappa>0$, with constants depending only on $\kappa$ and uniform spatial bounds, independently of the starting time and local temporal regularity norms. This includes the critical conductivity $a=x^2$ with noncompact thermal embedding. For autonomous powers $a=x^\alpha$ with the corresponding boundary realizations, exponential stability holds uniformly in $\alpha$ for $0<\alpha\le2$. For $\alpha>2$, concentrating approximate equilibria show that the semigroup norm equals $1$ at every time and preclude uniform decay even from the generator graph domain to the energy space. Proofs combine a common-domain construction, interior-average compensation and time averaging for Dirichlet conditions, and a strict Lyapunov functional for zero-flux conditions.
\end{abstract}
\maketitle

\section{Introduction}

Degenerate evolution equations arise when a principal coefficient
vanishes at a boundary or an interior point, affecting traces,
observability and stabilization.
For boundary-degenerate wave equations, Alabau-Boussouira,
Cannarsa and Leugering \cite{ABC2017} developed control and
stabilization results in weighted energy spaces.
Subsequent work treats linear and nonlinear stabilization in the
presence of drift terms \cite{FM2025nonlinear,FM2025linear},
simultaneous degeneracy of propagation and damping
\cite{AFIssa2025}, delayed feedbacks \cite{CFP2025,Liao2025},
and degenerate beams \cite{CF2024}, among others. Non-autonomous hyperbolic problems have also been studied through
boundary observability \cite{AFIobs} and boundary stabilization
in divergence and non-divergence form \cite{AFIdiv,AFInondiv}.
These results illustrate how the degeneracy, boundary conditions
and temporal variation of the coefficients influence the
propagation and dissipation of energy.

For thermoelastic systems, heat conduction dissipates mechanical
energy through its coupling with elasticity. Classical energy and
multiplier methods for this mechanism are developed in \cite{Kim}.
A related weakly degenerate system with a non-divergence thermal
operator and a time-varying boundary delay is considered in \cite{Salhi}.
Boundary degeneracy requires a weighted treatment of thermal traces
and of the transfer of dissipation to the mechanical component.
For time-dependent conductivities, this transfer must be quantified
uniformly along the evolution.

The non-autonomous parabolic framework most relevant here is
developed in three works by Akil, Fragnelli and Ismail
\cite{AFIRobin2025,AFIpar,AFIdivpar2026}.
The first treats separable conductivities of the form
$b(t)\widetilde a(x)$ in divergence and non-divergence form,
with particular emphasis on Robin conditions in the weakly
degenerate case \cite{AFIRobin2025}.
The subsequent papers consider general, nonseparable space--time
dependent coefficients in non-divergence and divergence form,
respectively, with Dirichlet or Robin boundary realizations
\cite{AFIpar,AFIdivpar2026}.
These works use weighted spaces and Kato's evolution theory for
well-posedness, and Carleman estimates for observability and
null-controllability. In particular, comparison with a fixed
spatial profile and a common-domain construction are already
available in this scalar parabolic setting.
Our well-posedness argument adapts this framework to the coupled
thermoelastic operator. The stabilization problem additionally
requires control of the conservative mechanical component by
thermal dissipation on the whole time half-line. The cited
finite-time controllability results do not directly supply
this estimate.

We consider the thermoelastic system
\begin{equation}\label{eq1}
\begin{cases}
u_{tt}-u_{xx}+\kappa\theta_x=0,
 &(x,t)\in(0,1)\times(0,\infty),\\
\theta_t-(a(x,t)\theta_x)_x+\kappa u_{xt}=0,
 &(x,t)\in(0,1)\times(0,\infty),\\
u(0,t)=u(1,t)=\theta(1,t)=0,
 &t>0,\\
u(x,0)=u_0(x),\quad u_t(x,0)=u_1(x),\quad
\theta(x,0)=\theta_0(x),
 &x\in(0,1),
\end{cases}
\end{equation}
where $\kappa>0$, $a(0,t)=0$ and $a(x,t)>0$ for $x>0$. Here $u$ is the displacement and $\theta$ is the temperature variation. At the degenerate endpoint, we impose
\begin{equation}\label{eq2}
\begin{cases}
\theta(0,t)=0,
 &\text{Dirichlet boundary condition},\\
(a\theta_x)(0,t)=0,
 &\text{zero-flux boundary condition}.
\end{cases}
\end{equation}
The boundary realization is fixed throughout each problem.
The Dirichlet condition is associated with integrability of the
reciprocal initial conductivity, and the zero-flux condition with
its divergence.

Writing $\widetilde a(x)=a(x,0)$, we assume uniform comparison
for well-posedness and an additional quadratic lower bound for decay,
in the form
\[
m\widetilde a(x)\le a(x,t)\le M\widetilde a(x),
\quad
\sup_{0<x\le1}\frac{x^2}{\widetilde a(x)}<\infty,
\]
together with the spatial and local temporal regularity specified
in Hypotheses~\ref{hyp2.2} and \ref{hyp2.3}.
The energy and thermal dissipation are
\[
E(t)=\frac12\int_0^1
\bigl(u_t^2+u_x^2+\theta^2\bigr)\dd x,
\quad
D(t)=\int_0^1a(x,t)|\theta_x|^2\dd x.
\]
For regular solutions, $E'(t)=-D(t)$. 
The coupling terms cancel in the energy balance,
so the thermal diffusion is the only source of dissipation. 
Since the conductivity does not enter the stored energy,
this identity contains no time derivative of $a$.
The main issue is to recover uniform control of the mechanical
energy from the thermal dissipation.

The closest autonomous thermoelastic predecessor is the rod
studied by Ammari, Hassine and Robbiano \cite{AHR}.
Under their assumptions on the conductivity and coupling,
they prove well-posedness and exponential stability for a fixed
conductivity $a=a(x)$ by semigroup and frequency-domain methods.
Their Dirichlet realization assumes $\frac1a\in L^1(0,1)$.
Their strongly degenerate zero-flux realization assumes
$a\in C^1[0,1]$, $\frac1a\notin L^1(0,1)$ and
$\frac1{\sqrt a}\in L^1(0,1)$, and their exponential-stability
theorem also requires $W\in L^1(0,1)$, where
$W(x)=\int_x^1\frac{\dd y}{a(y)}$ in that case
\cite[Theorem~3.3]{AHR}.
For $a=x^2$, both $\frac1{\sqrt a}=\frac1x$ and
$W(x)=\frac1x-1$ are nonintegrable.
Against this background, the present work provides three
stabilization results. First, for genuinely time-dependent
conductivities, it obtains exponential decay uniformly over
coefficient families with common spatial bounds, without
dependence on their local time-regularity norms.
Second, it includes the critical conductivity $x^2$, using
coercivity without a compact thermal embedding.
Third, for autonomous powers $x^\alpha$ with $\alpha>2$,
concentrating approximate equilibria exclude every uniform
decay factor tending to zero and give a quantitative time
scale for energy persistence.
These extensions do not assert that our sufficient weight
assumptions contain every coefficient class treated in \cite{AHR}.

Three difficulties arise in passing from these scalar and
autonomous results to the present stabilization problem.
First, for genuinely time-dependent conductivities there is
no single autonomous generator to which a frequency-domain
stability criterion can be applied directly.
The coupled domain must also account for the stress
$u_x-\kappa\theta$ as a whole, without imposing separate
$L^2$ regularity on $u_{xx}$ and $\theta_x$.
Multiplication estimates for $\frac{a}{\widetilde a}$ and
$\frac{\widetilde a}{a}$ give domain invariance and the operator
regularity needed by Kato's theorem. 
The construction uses the full thermal norm and does not require
the quadratic lower bound. That additional spatial condition
enters only the multiplier estimates for decay.

Second, in the Dirichlet case, the temperature primitive
produces an endpoint heat-flux term.
The vanishing temperature trace does not eliminate this flux,
and its trace cannot be bounded by the dissipation alone.
An interior-average correction removes the endpoint flux,
while time averaging controls the resulting low-order
mechanical terms.
Averaging the actual flux $a\theta_x$, rather than separating
its factors, avoids differentiating the conductivity in time.
In the zero-flux case, the left endpoint term vanishes directly,
and a mechanical multiplier compensates the remaining boundary
contribution to yield a strict Lyapunov inequality.

Third, at the critical coefficient $a=x^2$, the thermal
embedding into $L^2(0,1)$ is not compact.
An argument relying on this compactness is therefore unavailable.
We instead use a weighted Hardy estimate and direct multiplier
bounds to recover the coercivity needed for decay.
Beyond this critical exponent, concentrating states with
constant stress make the diffusion residual arbitrarily small,
providing the obstruction to uniform decay.

The main conclusions are as follows.
\begin{enumerate}
\item Under Hypotheses~\ref{hyp2.2} and \ref{hyp2.3}, the system
generates a contraction evolution family on
\[
\mathcal H=H_0^1(0,1)\times L^2(0,1)\times L^2(0,1)
\]
and is globally well posed for energy data. 
Under the additional Hypothesis \ref{hybound},
for every fixed $\kappa>0$ there exist $C\ge1$ and $\omega>0$
such that
\[
E(t)\le Ce^{-\omega(t-s)}E(s),
\quad t\ge s\ge0.
\]
The constants depend only on $\kappa$ and the uniform spatial
bounds, independently of the starting time and the local
time-regularity norms in Hypothesis~\ref{hyp2.3}.
In particular, no global bound on these time-regularity norms
is required.

\item For $a=x^\alpha$, with the Dirichlet condition at zero
when $0<\alpha<1$ and the natural zero-flux condition when
$\alpha\ge1$, exponential stability holds for $0<\alpha\le2$.
For fixed $\kappa>0$, the decay constants can be chosen
independently of $\alpha$ throughout this range.
At $\alpha=2$, the thermal embedding is not compact, but
the Hardy estimate
\[
\|z\|^2\le4\|xz_x\|^2,\quad z\in V_2,
\]
provides the coercivity needed for decay.

\item For $\alpha>2$, the autonomous operator still generates
a contraction semigroup $S_\alpha(t)$.
Concentrating states with constant stress and vanishing
diffusion residual yield
\[
\|S_\alpha(t)\|_{\mathcal L(\mathcal H)}=1,
\quad t\ge0.
\]
Thus no uniform decay factor tending to zero exists.
The same obstruction persists when the initial data are
controlled by the generator graph norm.
Together, these results identify the sharp threshold
$\alpha=2$ for uniform decay in the autonomous power family
with the stated boundary realizations.
\end{enumerate}

Section~\ref{sec:prelim} states the assumptions, functional setting
and main results.
Section~\ref{sec:proofs} establishes the weighted estimates,
the common operator domain and global well-posedness.
Section~\ref{sec:decay} proves exponential decay for both boundary
conditions.
Section~\ref{sec:powers} establishes the power threshold and
examines the critical noncompact case.
Section~\ref{sec:extensions} gives examples and discusses
the scope of the assumptions.

\section{Assumptions, functional setting and main results}\label{sec:prelim}

\subsection{Classical degeneracy and the coefficient assumptions}
We use the following uniform degeneracy conditions, with an
absolute derivative bound. For related terminology, see
\cite[Definition~2.3]{AFIpar}.
\begin{definition}[Classical uniform weak and strong degeneracy]
\label{def2.1}
Let $a(\cdot,t)\in C[0,1]\cap C^1((0,1])$ for every $t\ge0$,
with $a(0,t)=0$ and $a(x,t)>0$ for $0<x\le1$. Define
\begin{equation*}
K=\sup_{t\ge0}\sup_{0<x\le1}
  \frac{|x a_x(x,t)|}{a(x,t)}.
\end{equation*}
\begin{enumerate}
\item We call $a$ classically weakly degenerate, denoted
\textup{(WD)}, if $K<1$.
\item We call $a$ classically strongly degenerate, denoted
\textup{(SD)}, if $1\le K<2$ and
$a(\cdot,t)\in C^1[0,1]$ for every $t\ge0$.
\end{enumerate}
\end{definition}

The classical cases imply the corresponding weight conditions below:
\textup{(WD)} gives the Dirichlet realization, whereas \textup{(SD)}
gives the zero-flux realization. Remark~\ref{rem2.4} proves
these implications. The weight assumptions also admit coefficients
outside the classical classes.

\begin{hypothesis}
\label{hyp2.2}
For every $t\ge0$, assume
\[
a(\cdot,t)\in C[0,1]\cap W^{1,\infty}_{\rm loc}((0,1]),
\quad a(0,t)=0,\quad a(x,t)>0\quad\text{for }x>0.
\]
Set $\widetilde a(x)=a(x,0)$. There are constants
$0<m\le1\le M<\infty$ such that
\begin{equation}\label{eq3}
m\widetilde a(x)\le a(x,t)\le M\widetilde a(x),
\quad\text{for }0\le x\le1,\quad t\ge0.
\end{equation}
Choose one of the following two cases, with its fixed boundary realization:
\begin{enumerate}
\item \textsf{Dirichlet boundary realization.} The initial conductivity satisfies
\begin{equation}\label{eq4}
\int_0^1\frac{\dd x}{\widetilde a(x)}<\infty.
\end{equation}
\item \textsf{Zero-flux boundary realization.} The initial conductivity satisfies
\begin{equation}\label{eq5}
\int_0^1\frac{\dd x}{\widetilde a(x)}=\infty.
\end{equation}
\end{enumerate}
\end{hypothesis}

Under Hypothesis~\ref{hyp2.2}, the constant
$\widetilde{M}=M\norm{\widetilde a}_\infty$ is finite and independent of time.

\begin{hypothesis}
\label{hyp2.3}
For each finite $T>0$, assume the following temporal regularity
in the corresponding boundary case:
\begin{enumerate}
\item \textsf{Dirichlet boundary case.} 
\begin{equation}\label{eq7}
\frac{a}{\widetilde a}
 \in W^{1,\infty}(0,T;H^1(0,1)).
\end{equation}
\item \textsf{Zero-flux boundary case.} 
\begin{equation}\label{eq8}
\frac{a}{\widetilde a},\
 x\partial_x\!\left(\frac{a}{\widetilde a}\right)
 \in W^{1,\infty}(0,T;L^\infty(0,1)).
\end{equation}
\end{enumerate}
\end{hypothesis}
Throughout, we use the continuous time representatives of the
coefficient ratios supplied by Hypothesis~\ref{hyp2.3}.
Their initial time traces are required to equal $1$ for
$\frac{a}{\widetilde a}$ and, in the zero-flux case, $0$ for
$x\partial_x(\frac{a}{\widetilde a})$.
All pointwise-in-time statements refer to these representatives.
The norms in these regularity assumptions may depend on $T$;
no uniform bound as $T\to\infty$ is required.

For the decay estimates, we additionally impose 
\begin{hypothesis}\label{hybound}
The following quadratic weight bound holds.
\begin{equation*}
\sup_{0<x\le1}\frac{x^2}{\widetilde a(x)}<\infty.
\end{equation*}
\end{hypothesis}
When this condition holds, we use the constants
\begin{equation}\label{eq6}
\begin{aligned}
M_1&=\frac1m\sup_{0<x\le1}\frac{x^2}{\widetilde a(x)},
&M_2&=4M_1.
\end{aligned}
\end{equation}
Both constants are finite, positive and independent of time.

\begin{remark}
\label{rem2.4}
\begin{enumerate}
\item If
$m_*a_*(x)\le a(x,t)\le M_*a_*(x)$
for a fixed positive function $a_*$ on $(0,1]$ and constants
$0<m_*\le M_*<\infty$, then comparison at $t=0$ yields
\begin{equation*}
\frac{m_*}{M_*}a(x,0)\le a(x,t)
\le \frac{M_*}{m_*}a(x,0).
\end{equation*}
Hence we may take $\widetilde a(x)=a(x,0)$ without losing
the two-sided comparison bounds. 
\item For a classical coefficient with finite index $K$, integration of
$|\widetilde a'|/\widetilde a\le K/x$ on $[x,1]$ yields
$
\log\frac{\widetilde a(1)}{\widetilde a(x)}
\le K\log\frac1x,
$ 
further, 
$
\widetilde a(x)\ge\widetilde a(1)x^K.
$ 
Consequently, $K<1$ implies \eqref{eq4} and
Hypothesis \ref{hybound}, with
\begin{equation}\label{eq9}
\int_0^1\frac{\dd x}{\widetilde a(x)}
\le\frac1{\widetilde a(1)(1-K)},
\quad
\sup_{0<x\le1}\frac{x^2}{\widetilde a(x)}
\le\frac1{\widetilde a(1)}.
\end{equation}
In the classically strongly degenerate case, $K<2$ still gives
the estimate for $\sup_{0<x\le1}x^2/\widetilde a(x)$ in
\eqref{eq9}.
The additional assumptions $\widetilde a\in C^1[0,1]$ and
$\widetilde a(0)=0$ imply
$
0<\widetilde a(x)\le\|\widetilde a'\|_\infty x,
$
which gives 
$
\int_0^1\frac{\dd x}{\widetilde a(x)}=\infty.
$
Here $\|\widetilde a'\|_\infty>0$ because
$\widetilde a(0)=0$ and $\widetilde a(1)>0$. 
Thus \textup{(SD)} implies \eqref{eq5} and
Hypothesis \ref{hybound}. 
The divergence of the reciprocal integral follows from
$C^1$ regularity at zero and $\widetilde a(0)=0$.
More generally, any coefficient $\widetilde a\in C^1[0,1]$
that vanishes only at zero and satisfies
$\widetilde a(x)\ge c x^2$ on $(0,1]$ for some $c>0$ 
satisfies \eqref{eq5} and Hypothesis \ref{hybound}.

\item Hypotheses~\ref{hyp2.2} and \ref{hyp2.3} imply that, for every
finite $T>0$,
\begin{equation}\label{eq10}
\left\|\frac{a_t}{a}\right\|_{L^\infty(0,T;L^\infty(0,1))}
\le \frac1m
\left\|\frac{a_t}{\widetilde a}\right\|_{L^\infty(0,T;L^\infty(0,1))}
<\infty.
\end{equation}
In the Dirichlet case, this also uses the embedding
$H^1(0,1)\hookrightarrow L^\infty(0,1)$. 
\end{enumerate}
\end{remark}

\subsection{Energy space and solutions}

We use real-valued functions. Unless a space is specified, $\norm{\cdot}$
denotes the norm in $L^2(0,1)$, and spatial integrals are over $(0,1)$.

Define the fixed thermal space
\begin{equation*}
V=\begin{cases}
\left\{\begin{aligned}
z\in L^2(0,1):\ z\in AC([0,1]), \sqrt{\widetilde a}\,z_x\in L^2,\ z(0)=z(1)=0
\end{aligned}\right\}, \text{Dirichlet case},\\[2mm]
\left\{\begin{aligned}
z\in L^2(0,1):\ z\in AC_{\rm loc}((0,1]), \sqrt{\widetilde a}\,z_x\in L^2,\ z(1)=0
\end{aligned}\right\}, \text{zero-flux case},
\end{cases}
\end{equation*}
with
$\norm{z}_V^2=\norm{z}^2+\int_0^1\widetilde a z_x^2\dd x$.
Indeed, comparison results in the explicit norm equivalence
\[
\min\{1,m\}\norm z_V^2
\le\norm z^2+\int_0^1a(x,t)z_x^2\dd x
\le\max\{1,M\}\norm z_V^2.
\]
The space $V$ is Hilbert and is continuously and densely
embedded in $L^2(0,1)$.
Completeness follows from local $H^1$ convergence and
convergence of the weighted derivatives; in the Dirichlet
case, $1/\widetilde a\in L^1(0,1)$ also ensures
preservation of the trace at zero. 

We now define the energy space $\mathcal H$ by
\[\mathcal H=H_0^1(0,1)\times L^2(0,1)\times L^2(0,1)\]
equipped with the norm
\begin{equation}\label{eq11}
\norm{(u,v,\theta)}_{\mathcal H}^2=\norm{u_x}^2+\norm v^2+\norm\theta^2.
\end{equation}

Denote $v=u_t$ and $U=(u,v,\theta)$, the system can be written as  $U_t=\mathcal A(t)U$, where
\begin{equation}\label{eq12}
\mathcal A(t)U=\bigl(v,(u_x-\kappa\theta)_x,(a(\cdot,t)\theta_x)_x-\kappa v_x\bigr).
\end{equation}
Its domain $\mathcal D(t)$ consists of the triples in $\mathcal H$ such that
\begin{equation}\label{eq13}
u,v\in H_0^1(0,1),\quad \theta\in V,\quad
u_x-\kappa\theta\in H^1(0,1),\quad
a(\cdot,t)\theta_x\in H^1(0,1),
\end{equation}
with the additional condition $(a(\cdot,t)\theta_x)(0)=0$ {in the zero-flux case}.
The stress condition in \eqref{eq13} is imposed on
$u_x-\kappa\theta$ as a whole; we do not separately assume
$u_{xx},\theta_x\in L^2(0,1)$. We prove below that
$\mathcal D(t)=\mathcal D(0)$ and write this common domain as $\mathcal D$,
equipped with the graph norm of $\mathcal A(0)$.

\begin{definition}[Solutions]\label{def2.5}
Let $s\ge0$. A strong solution with datum $U_s\in\mathcal D$ is a function
\[
U\in C([s,T];\mathcal D)\cap C^1([s,T];\mathcal H)
\quad\text{for every }T>s
\]
such that $U_t=\mathcal A(t)U$ and $U(s)=U_s$.
For $U_s\in\mathcal H$, an energy solution is the limit in
$C([s,T];\mathcal H)$, for every $T>s$, of strong solutions whose data in
$\mathcal D$ converge to $U_s$ in $\mathcal H$. 
\end{definition}

\subsection{Main results}

\begin{theorem}[Global well-posedness]\label{thm2.6}
Assume Hypotheses~\ref{hyp2.2} and \ref{hyp2.3}, fix
$\kappa>0$, and impose the corresponding thermal boundary
condition in \eqref{eq2}. The family \eqref{eq12} has a dense
common domain $\mathcal D$ 
and generates a unique contraction evolution family on $\mathcal H$. For every $s\ge0$ and $U_s\in\mathcal H$, problem
\eqref{eq1}--\eqref{eq2}, initialized at time $s$, has a
unique global energy solution. It satisfies
\begin{equation}\label{eq14}
E(t)+\int_s^tD(\tau)\dd\tau=E(s),\quad t\ge s,
\end{equation}
and $\theta\in L^2(s,T;V)$ for every finite $T>s$. If
$U_s\in\mathcal D$, the solution is strong in the sense of
Definition~\ref{def2.5}. For any two solutions,
\begin{equation}\label{eq15}
\norm{U(t)-\widehat U(t)}_{\mathcal H}
\le\norm{U(s)-\widehat U(s)}_{\mathcal H},\quad t\ge s.
\end{equation}
\end{theorem}

\begin{theorem}[Uniform energy decay]\label{thm2.7}
Assume Hypotheses~\ref{hyp2.2}, \ref{hyp2.3} and \ref{hybound},
fix $\kappa>0$, and impose the corresponding boundary condition
\eqref{eq2}. There are constants $C\ge1$ and $\omega>0$,
depending only on $\kappa,M_1,M_2,\widetilde M$ and independent
of the starting time and the local time-regularity norms,
such that every energy solution satisfies
\begin{equation}\label{eq16}
E(t)\le C e^{-\omega(t-s)}E(s),\quad t\ge s\ge0.
\end{equation}
\end{theorem}

The following corollary follows directly from
Remark~\ref{rem2.4} and
Theorems~\ref{thm2.6} and \ref{thm2.7}.

\begin{corollary}[Classical weak and strong degeneracy]
\label{cor2.8}
Let $\kappa>0$. Suppose that $a$ is classically \textup{(WD)}
or \textup{(SD)} and satisfies \eqref{eq3} with
$\widetilde a=a(\cdot,0)$.
Under the Dirichlet boundary condition in the \textup{(WD)} case
or the zero-flux boundary condition in the \textup{(SD)} case,
and the corresponding time regularity in Hypothesis~\ref{hyp2.3},
problem 
\eqref{eq1}--\eqref{eq2} is globally well posed in the sense of
Theorem~\ref{thm2.6}, and every energy solution satisfies
the exponential decay estimate \eqref{eq16}.
\end{corollary}

\begin{theorem}[Sharp threshold for autonomous powers]
\label{thm2.9}
Let $\alpha>0$ and $\kappa>0$, and set $a=x^\alpha$.
Impose $\theta(0,t)=0$ for $0<\alpha<1$ and
$(x^\alpha\theta_x)(0,t)=0$ for $\alpha\ge1$,
in addition to $u(0,t)=u(1,t)=\theta(1,t)=0$.
The operator $(\mathcal A_\alpha,\mathcal D_\alpha)$ defined in
\eqref{119}--\eqref{eq62} generates a
contraction semigroup $S_\alpha(t)$ on $\mathcal H$.
\begin{enumerate}
\item If $0<\alpha\le2$, there are $C\ge1$ and $\omega>0$ such that
$E(t)\le Ce^{-\omega(t-s)}E(s)$ for all energy solutions and $t\ge s$.
For fixed $\kappa$, these constants can be chosen independently
of $\alpha$ in this range.
\item If $\alpha>2$, then
\begin{equation}\label{eq17}
\|S_\alpha(t)\|_{\mathcal L(\mathcal H)}=1
\quad\text{for every }t\ge0.
\end{equation}
In particular, no function $r(t)\to0$ can satisfy
$E(t)\le r(t)E(0)$ for all initial data in $\mathcal H$.
\end{enumerate}
\end{theorem}

\section{Analytic tools and global well-posedness}\label{sec:proofs}
Under Hypothesis \ref{hybound}, the constants in
\eqref{eq6} give
\begin{equation}\label{eq18}
a(x,t)\ge m\widetilde a(x)\ge\frac{x^2}{M_1},
\quad \text{for }0<x\le1.
\end{equation}
The upper bound $a\le \widetilde{M}$ holds under
Hypothesis~\ref{hyp2.2} alone. The same hypothesis gives
the following distinction between the two boundary realizations:
\begin{equation}\label{eq19}
\begin{aligned}
\text{Dirichlet case}:\quad
\int_0^1\frac{1}{a(x,t)}\dd x
&\le\frac1m\int_0^1\frac{1}{\widetilde a}\dd x<\infty,\\
\text{Zero-flux case}:\quad
\int_0^1\frac{1}{a(x,t)}\dd x
&\ge\frac1M\int_0^1\frac{1}{\widetilde a}\dd x=\infty.
\end{aligned}
\end{equation}

\begin{lemma}[Weighted estimates and endpoint behavior]\label{lem3.1}
Under Hypotheses~\ref{hyp2.2} and \ref{hybound},
for every $z\in V$ and $t\ge0$,
\begin{equation}\label{eq20}
\|xz_x\|^2\le M_1\int_0^1az_x^2\dd x,\quad
\|z\|^2\le M_2\int_0^1az_x^2\dd x,\quad
\|az_x\|^2\le \widetilde{M}\int_0^1az_x^2\dd x\end{equation}
Moreover,
\begin{equation}\label{eq21}
\lim_{x\downarrow0}x|z(x)|^2=0.
\end{equation}
\end{lemma}

\begin{proof}
Fix $t\ge0$. The first estimate in \eqref{eq20} follows
from \eqref{eq18}, and the third follows from
$a\le\widetilde M$. To prove the second estimate, let
$z\in V$. 
For $0<\varepsilon<1$, local $H^1$ regularity, $z(1)=0$ and  Cauchy-Schwarz  yield
\[
\int_\varepsilon^1z^2\dd x
=-\varepsilon z(\varepsilon)^2-2\int_\varepsilon^1xzz_x\dd x
\le2\|z\|_{L^2(\varepsilon,1)}\|xz_x\|_{L^2(\varepsilon,1)}.
\]
Letting $\varepsilon\downarrow0$ yields
\begin{equation}
\label{118}
\|z\|^2\le4\|xz_x\|^2\le M_2\int_0^1a|z_x|^2\dd x.
\end{equation}
Moreover, $(xz^2)'=z^2+2z(xz_x)\in L^1(0,1)$, so $xz^2$ has a
finite limit $\iota$ at zero. If $\iota>0$, then $z(x)^2\ge \iota/(2x)$
for all sufficiently small $x>0$, contradicting
$z\in L^2(0,1)$ since $1/x$ is not integrable near zero. This proves \eqref{eq21}.
\end{proof}
Let us recall that the Hardy's inequality for $z\in H^1(0,1)$ with $z(0)=0$ is
\begin{equation}\label{eq25}
\left\|\frac zx\right\|^2\le4\|z_x\|^2.
\end{equation}
\begin{lemma}[Flux integration by parts]\label{lem3.2}
Assume Hypothesis~\ref{hyp2.2} and fix $t\ge0$.
If $\theta,z\in V$ and $a(\cdot,t)\theta_x\in H^1(0,1)$, then
\begin{equation}\label{eq22}
\int_0^1(a\theta_x)_xz\dd x
=-\int_0^1a(\cdot,t)\theta_xz_x\dd x.
\end{equation}
If, in addition, Hypothesis \ref{hybound} holds,
then every $U\in\mathcal D(t)$ satisfies
\begin{equation}\label{eq23}
\lim_{x\downarrow0}x|u_x(x)|^2=0.
\end{equation}
\end{lemma}

\begin{proof}
The boundary cancellation follows the argument of
\cite[Proposition~2.1]{AHR}. We establish the Green identity
without using Hypothesis \ref{hybound}.
In the Dirichlet case, the flux is bounded and $z(0)=0$,
so $(a\theta_x)(x)z(x)\to0$ as $x\downarrow0$.
In the zero-flux case,
\[
\int_0^1\frac{|a\theta_x|^2}{a}\dd x
=\int_0^1a|\theta_x|^2\dd x<\infty.
\]
If $(a\theta_x)(0)\ne 0$, by continuity there exists $\delta>0$ such that
$
|(a\theta_x)(x)|\ge \frac{|(a\theta_x)(0)|}{2}
$ for $0<x<\delta.$ 
Then by the reciprocal divergence in
\eqref{eq19}, 
\[
\int_0^\delta \frac{|a\theta_x|^2}{a}\,dx
\ge \frac{|(a\theta_x)(0)|^2}{4}\int_0^\delta \frac{dx}{a}
=\infty,
\]
which yields a contradiction. 
Moreover, $(a\theta_x)z\in L^1(0,1)$ and
\[
\partial_x\bigl((a\theta_x)z\bigr)
=(a\theta_x)_xz+(\sqrt a\,\theta_x)(\sqrt a\,z_x)
\in L^1(0,1)
\]
by Cauchy--Schwarz and comparison with $\widetilde a$.
Hence the product belongs to $W^{1,1}(0,1)$ and has a finite
limit at zero. Hardy's inequality \eqref{eq25}, applied to
the zero-trace flux, gives
\[
\int_0^1\frac{|(a\theta_x)z|}{x}\dd x
\le\left\|\frac{a\theta_x}{x}\right\|\|z\|
\le2\|(a\theta_x)_x\|\|z\|<\infty.
\]
Since the product has a finite limit at zero, a nonzero limit
would imply the existence of $c>0$ and $\delta_0>0$ such that
\[
|(a\theta_x)(x)z(x)|\ge c,
\quad 0<x<\delta_0.
\]
This would give
\[
\int_0^{\delta_0}\frac{|(a\theta_x)z|}{x}\dd x
\ge c\int_0^{\delta_0}\frac{\dd x}{x}
=\infty,
\]
contradicting the preceding integrability estimate.
Hence $(a\theta_x)(x)z(x)\to0$ as $x\downarrow0$
also in the zero-flux case. 
Integration by parts on $(\delta,1)$, using $z(1)=0$, yields
\[
\int_\delta^1(a\theta_x)_xz\dd x
=-(a\theta_x)(\delta)z(\delta)
-\int_\delta^1a\theta_xz_x\dd x.
\]
Letting $\delta\downarrow0$, the boundary term vanishes and
both integrals converge by their $L^1$ integrability,
which proves \eqref{eq22}.

For the mechanical trace, assume in addition
Hypothesis~\ref{hybound}. The embedding
$H^1(0,1)\hookrightarrow L^\infty(0,1)$, applied to
$u_x-\kappa\theta$, and \eqref{eq21} give
\[
x|u_x(x)|^2
\le2x\|u_x-\kappa\theta\|_\infty^2
  +2\kappa^2x|\theta(x)|^2\to0\quad\text{as }x\downarrow0,
\]
which proves \eqref{eq23}.
\end{proof}

\begin{lemma}
\label{lem3.3}
Let $z\in H^1(0,1)$ with $z(0)=0$.  Fix $t\ge0$ and suppose that $a/\widetilde a$ has a local weak
derivative, with
\[
0<m\le\frac{a}{\widetilde a}\le M,\quad
x\partial_x\!\left(\frac{a}{\widetilde a}\right)\in L^\infty(0,1).
\]
Then both $(a/\widetilde a)z$ and $(\widetilde a/a)z$
belong to $H^1(0,1)$ and have zero trace at $x=0$. Moreover,
\begin{equation}\label{eq24}
\left\|\frac{a}{\widetilde a}z\right\|_{H^1}
\le2\left(
\left\|\frac{a}{\widetilde a}\right\|_\infty
+\left\|x\partial_x\!\left(\frac{a}{\widetilde a}\right)\right\|_\infty
\right)\|z\|_{H^1}.
\end{equation}
The corresponding estimate for $(\widetilde a/a)z$ is obtained
by replacing every occurrence of $a/\widetilde a$ in
\eqref{eq24} by $\widetilde a/a$.\end{lemma}

\begin{proof}

%Indeed, since $z(0)=0$, by Cauchy--Schwarz, 
%\[|z(x)|^2=\left|\int_0^x y^{-1/4}y^{1/4}z_x(y)\dd y\right|^2\le2\sqrt{x}\int_0^x\sqrt{y}|z_x(y)|^2\dd y.\]
%Multiplying by $1/x^{2}$, integrating, and applying Tonelli'stheorem, we obtain
%\[\int_0^1\frac{|z(x)|^2}{x^2}\dd x\le2\int_0^1\sqrt{y}|z_x(y)|^2\int_y^1x^{-3/2}\dd x\dd y=4\int_0^1(1-\sqrt{y})|z_x(y)|^2\dd y\le4\|z_x\|^2.\] 

The product rule yields
$
\partial_x\!\left(\frac{a}{\widetilde a}z\right)
=\frac{a}{\widetilde a}z_x
+x\partial_x\!\left(\frac{a}{\widetilde a}\right)\frac{z}{x}.
$
Hence, by \eqref{eq25},
\[
\left\|\partial_x\!\left(\frac{a}{\widetilde a}z\right)\right\|
\le\left(
\left\|\frac{a}{\widetilde a}\right\|_\infty
+2\left\|x\partial_x\!\left(\frac{a}{\widetilde a}\right)\right\|_\infty
\right)\|z_x\|.
\]
Together with
$\|(a/\widetilde a)z\|\le\|a/\widetilde a\|_\infty\|z\|$,
this proves \eqref{eq24}.

The trace is zero since
$|(a/\widetilde a)(x)z(x)|\le M|z(x)|\to0$ as $x\downarrow0$.
The same argument applies to $\widetilde a/a$. Indeed, the
local weak chain rule gives
\[
x\partial_x\!\left(\frac{\widetilde a}{a}\right)
=-\left(\frac{\widetilde a}{a}\right)^2
x\partial_x\!\left(\frac a{\widetilde a}\right).
\]
Consequently,
\[
\left\|\frac{\widetilde a}{a}\right\|_\infty\le\frac1m,
\qquad
\left\|x\partial_x\!\left(\frac{\widetilde a}{a}\right)\right\|_\infty
\le\frac1{m^2}
\left\|x\partial_x\!\left(\frac a{\widetilde a}\right)\right\|_\infty.
\]
\end{proof}

We use the following form of Kato's evolution theorem
\cite[Theorem~1.2]{Kato}; see also \cite[Theorem~3]{AFIobs}
and \cite[Theorem~3.3]{AFIpar}.

\begin{theorem}\label{thm3.4}
Let $T>0$ be fixed, $\mathcal{H}$ a Hilbert space and $\mathcal{A}(t): D(\mathcal{A}(t))\subset \mathcal{H}\to \mathcal{H}$ a linear operator satisfying the following hypotheses:

\begin{itemize}
\item[$(H_1)$] $D(\mathcal{A}(t)) = D(\mathcal{A}(0)),\ \forall t\in [0,T]$;
\item[$(H_2)$] $D(\mathcal{A}(0))$ is a dense subset of $\mathcal{H}$;
\item[$(H_3)$] $\forall t\in [0,T]$, $\mathcal{A}(t)$ generates a strongly continuous semigroup on $\mathcal{H}$, and the family $\mathscr{A}=\{\mathcal{A}(t)\}_{t\in[0,T]}$ is stable with stability constants $C$ and $p$ independent of $t$

(i.e. the semigroup $(S_t(s))_{s\ge 0}$ generated by $\mathcal{A}(t)$ satisfies
\[
\left\|\prod_{j=1}^k S_{t_j}(s_j)u\right\|_{\mathcal{H}}
\le C e^{p\sum_{j=1}^k s_j}\|u\|_{\mathcal{H}},\quad
\text{for all }u\in\mathcal{H},\ s_j\ge 0,\ 0\le t_1\le \dots t_k\le T,\ k\in\mathbb{N});
\]
\item[$(H_4)$] $\frac{\partial}{\partial t}\mathcal{A}\in L^\infty(0,T;B(D(\mathcal{A}(0)),\mathcal{H}))$, which is the space of equivalent classes of essentially bounded, strongly measurable functions from $[0,T]$ into the set of bounded operators $B(D(\mathcal{A}(0)),\mathcal{H})$.
\end{itemize}

Then, the problem
\[
\begin{cases}
\dot u(t) = \mathcal{A}(t)u(t), & t\in (0,T),\\
u(0)=u_0\in \mathcal{H},
\end{cases}
\]
has a unique solution $u\in C([0,T],\mathcal{H})$ for all $u_0\in \mathcal{H}$. Moreover, if $u_0\in D(\mathcal{A}(0))$, then $u\in C([0,T];D(\mathcal{A}(0)))\cap C^1([0,T];\mathcal{H})$.
\end{theorem}

\begin{proposition}[Maximal dissipativity]
\label{prop3.5}
Under Hypothesis~\ref{hyp2.2}, for every $t\ge0$,
$\mathcal D(t)$ is dense in $\mathcal H$ and $\mathcal A(t)$ is
maximal dissipative. It generates a contraction semigroup $S_t(r)$
on $\mathcal H$, with
\[
\|S_t(r)\|_{\mathcal L(\mathcal H)}\le1,
\quad r\ge0.
\]
Moreover, for $U=(u,v,\theta)\in\mathcal D(t)$,
\begin{equation}\label{eq26}
\pair{\mathcal A(t)U}{U}_{\mathcal H}
=-\int_0^1a(x,t)\theta_x^2\dd x.
\end{equation}
\end{proposition}

\begin{proof}
Fix $t\ge0$. Since $a(\cdot,t)$ is locally Lipschitz,
$C_c^\infty(0,1)^3\subset\mathcal D(t)$, so $\mathcal D(t)$
is dense in $\mathcal H$. 

For $U=(u,v,\theta)\in\mathcal D(t)$,
integration by parts, using $v\in H_0^1(0,1)$ and
\eqref{eq22} with $z=\theta$, shows
\[
\begin{aligned}
\langle\mathcal A(t)U,U\rangle_{\mathcal H}
&=\int_0^1u_xv_x\dd x
-\int_0^1(u_x-\kappa\theta)v_x\dd x
-\int_0^1a|\theta_x|^2\dd x
-\kappa\int_0^1v_x\theta\dd x\\
&=-\int_0^1a|\theta_x|^2\dd x\le0,
\end{aligned}
\]
which proves \eqref{eq26}.

For $(f,g,h)\in\mathcal H$, set $v=u-f$.
The equation $(I-\mathcal A(t))U=(f,g,h)$ reduces to
\begin{equation*}
\begin{cases}
u-(u_x-\kappa\theta)_x=f+g,\\
\theta-(a\theta_x)_x+\kappa u_x=h+\kappa f_x.
\end{cases}
\end{equation*}
On $X=H_0^1(0,1)\times V$, use
\[
\|(u,\theta)\|_X^2
=\|u\|^2+\|u_x\|^2+\|\theta\|^2
 +\|\sqrt{\widetilde a}\,\theta_x\|^2
\]
and the continuous bilinear form
\begin{equation}\label{eq27}
\begin{aligned}
B_t((u,\theta),(\varphi,\psi))
=\int_0^1(u\varphi+u_x\varphi_x-\kappa\theta\varphi_x)\dd x+\int_0^1(\theta\psi+a\theta_x\psi_x+\kappa u_x\psi)\dd x.
\end{aligned}
\end{equation}
The opposite coupling signs result in, for every $\kappa>0$,
\begin{equation*}
\begin{aligned}
B_t((u,\theta),(u,\theta))
=\|u\|^2+\|u_x\|^2+\|\theta\|^2
  +\int_0^1a\theta_x^2\dd x 
  \ge\min\{1,m\}\|(u,\theta)\|_X^2.
\end{aligned}
\end{equation*}
This coercivity uses the full norm of $V$. 
The functional
\[
(\varphi,\psi)\longmapsto
\int_0^1\bigl((f+g)\varphi+(h+\kappa f_x)\psi\bigr)\dd x
\]
is continuous on $X$. Thus Lax--Milgram
\cite[Corollary~5.8]{Brezis} yields a unique solution $
(u,\theta)\in H_0^1(0,1)\times V$. Further, $v\in H_0^1(0,1)$. 
Compactly supported tests imply
\[
(u_x-\kappa\theta)_x=v-g\in L^2(0,1),\quad
(a\theta_x)_x=\theta+\kappa v_x-h\in L^2(0,1).
\]
The stress is in $L^2(0,1)$ and
$\|a\theta_x\|^2\le \widetilde{M}\int_0^1a|\theta_x|^2\dd x$,
so both stress and flux belong to $H^1(0,1)$.
The reciprocal-divergence argument in Lemma~\ref{lem3.2}
gives the zero-flux trace when required.
Consequently, $U\in\mathcal D(t)$ and
$I-\mathcal A(t)$ is onto. Together with dissipativity, this proves
maximal dissipativity and closedness. The Lumer--Phillips theorem
\cite[Theorem~II.3.15]{EN} gives the contraction semigroup.
\end{proof}

We equip $\mathcal D(t)$ with the graph norm
\begin{equation*}
\|U\|_{\mathcal D(t)}
=\|U\|_{\mathcal H}+\|\mathcal A(t)U\|_{\mathcal H}.
\end{equation*}
Since $\mathcal A(t)$ is closed, this makes $\mathcal D(t)$ a
Banach space. In applying Theorem~\ref{thm3.4}, we use this norm
at $t=0$ on $D(\mathcal A(0))=\mathcal D(0)$.
The notation $B(\mathcal D(0),\mathcal H)$ in that theorem means
the same bounded-operator space as $\mathcal L(\mathcal D(0),\mathcal H)$.

\begin{proposition}\label{prop3.6}
Under Hypotheses~\ref{hyp2.2} and \ref{hyp2.3}, for $t\ge0$, 
$
\mathcal D(t)=\mathcal D(0)=:\mathcal D.
$
Write $\|U\|_{\mathcal D}=\|U\|_{\mathcal D(0)}$.
For every $T>0$, there exists a constant $C_T\ge1$,
depending only on $T$, $\kappa$, $\widetilde a$, the comparison
constants $m,M$, and the time-regularity norms on $(0,T)$
in Hypothesis~\ref{hyp2.3}, such that
\begin{equation}\label{eq28}
C_T^{-1}\|U\|_{\mathcal D}
\le\|U\|_{\mathcal D(t)}
\le C_T\|U\|_{\mathcal D},
\quad U\in\mathcal D,\quad t\in[0,T].
\end{equation}
Moreover,
\begin{equation*}
\mathcal A\in W^{1,\infty}(0,T;\mathcal L(\mathcal D,\mathcal H)).
\end{equation*}
Condition \textup{(H$_3$)} of Theorem~\ref{thm3.4} holds
with $C=1$ and $p=0$.
\end{proposition}

\begin{proof}
In the Dirichlet case, \eqref{eq7} and
$a/\widetilde a\ge m$ imply that both $a(\cdot,t)/\widetilde a$
and its reciprocal are uniformly bounded in $H^1(0,1)$
for $t\in[0,T]$. The one-dimensional $H^1$ product estimate
therefore yields
\begin{equation}
\label{eq29}
\left\|\frac{a(\cdot,t)}{\widetilde a}z\right\|_{H^1}
+\left\|\frac{\widetilde a}{a(\cdot,t)}z\right\|_{H^1}
\le C_T\|z\|_{H^1}.
\end{equation}
In the zero-flux case, the same estimate holds for
$z\in H^1(0,1)$ with $z(0)=0$ by
Lemma~\ref{lem3.3} and \eqref{eq8},
and both products retain zero trace.
Applying these estimates to the identities
\[
a(\cdot,t)\theta_x
=\frac{a(\cdot,t)}{\widetilde a}\,\widetilde a\theta_x,
\quad
\widetilde a\theta_x
=\frac{\widetilde a}{a(\cdot,t)}\,a(\cdot,t)\theta_x
\]
shows that 
\[
a(\cdot,t)\theta_x\in H^1(0,1)
\quad\Longleftrightarrow\quad
\widetilde a\theta_x\in H^1(0,1).
\]
All other domain conditions are independent of $t$, proving
$\mathcal D(t)=\mathcal D(0)=:\mathcal D$.

For $U=(u,v,\theta)\in\mathcal D$, dissipativity
\eqref{eq26} and the upper bound $a\le \widetilde{M}$ give
\[\|a\theta_x\|^2
\le \widetilde{M}\int_0^1a\theta_x^2\dd x
=-\widetilde{M}\langle\mathcal A(t)U,U\rangle_{\mathcal H}
\le \widetilde{M}\|\mathcal A(t)U\|_{\mathcal H}\|U\|_{\mathcal H},\]
further,
$\|a\theta_x\|\le \frac{\sqrt{\widetilde{M}}}{2}
 \bigl(\|U\|_{\mathcal H}+\|\mathcal A(t)U\|_{\mathcal H}\bigr).$ 
Recall \eqref{eq11} and \eqref{eq12}, one gets
\[\|(a\theta_x)_x\|
\le \|(a\theta_x)_x-\kappa v_x\|+\kappa\|v_x\|\le (1+\kappa)\|\mathcal A(t)U\|_{\mathcal H}.
\]
Applying \eqref{eq29} with
$z=a(\cdot,t)\theta_x$, we obtain 
\begin{equation}\label{eq30}
\|\widetilde a\theta_x\|_{H^1}
=\left\|\frac{\widetilde a}{a(\cdot,t)}
          \bigl(a(\cdot,t)\theta_x\bigr)\right\|_{H^1}
\le C_T\|a(\cdot,t)\theta_x\|_{H^1}
\le C_T\bigl(\|U\|_{\mathcal H}
+\|\mathcal A(t)U\|_{\mathcal H}\bigr).
\end{equation}
By the definition in \eqref{eq12}, 
$
\|(\mathcal A(t)-\mathcal A(s))U\|_{\mathcal H}
=
\left\|\partial_x\!\left(
[a(\cdot,t)-a(\cdot,s)]\theta_x\right)\right\|.
$ 
In the Dirichlet case, the $H^1(0,1)$ product estimate directly shows 
\begin{equation}
\label{eq31}
\|(\mathcal A(t)-\mathcal A(s))U\|_{\mathcal H}
\le
C\|\widetilde a\theta_x\|_{H^1}
\left\|
\frac{a(\cdot,t)-a(\cdot,s)}{\widetilde a}
\right\|_{H^1}.
\end{equation}
In the zero-flux case, the product estimate proved in
Lemma~\ref{lem3.3} shows
\begin{equation}
\label{eq32}
\begin{aligned}
\|(\mathcal A(t)-\mathcal A(s))U\|_{\mathcal H}
&\le
C\|\widetilde a\theta_x\|_{H^1}
\left\|
\frac{a(\cdot,t)-a(\cdot,s)}{\widetilde a}
\right\|_\infty\\
&\quad+
C\|\widetilde a\theta_x\|_{H^1}
\left\|x\partial_x\!\left(
\frac{a(\cdot,t)-a(\cdot,s)}{\widetilde a}
\right)\right\|_\infty.
\end{aligned}
\end{equation}
Combining
\eqref{eq31}, \eqref{eq32}  with \eqref{eq30} at time $s$
and using the triangle inequality yields
\[
\begin{aligned}
\|U\|_{\mathcal H}+\|\mathcal A(t)U\|_{\mathcal H}
\le \|U\|_{\mathcal H}+\|\mathcal A(s)U\|_{\mathcal H}
  +\|(\mathcal A(t)-\mathcal A(s))U\|_{\mathcal H} 
  \le C_T\bigl(
\|U\|_{\mathcal H}+\|\mathcal A(s)U\|_{\mathcal H}
\bigr).
\end{aligned}
\]
Interchanging $s,t$ and then taking $s=0$ proves
\eqref{eq28}. In particular,
\begin{equation}
\label{eq33}
\|\widetilde a\theta_x\|_{H^1}
\le C_T\|U\|_{\mathcal D}.
\end{equation}

We next verify the time regularity required in \textup{(H$_4$)}.  For fixed $U=(u,v,\theta)\in \mathcal D$, 
in the Dirichlet case, we have
\[
(\mathcal A(t)-\mathcal A(0))U
=
\left(0,0,\partial_x\!\left[
\left(\frac{a(\cdot,t)}{\widetilde a}-1\right)
\widetilde a\theta_x
\right]\right).
\]
The one-dimensional $H^1$ product estimate gives
\[
\left\|
\partial_x\!\left[
\left(\frac{a(\cdot,t)}{\widetilde a}-1\right)
\widetilde a\theta_x
\right]\right\|
\le
C\left\|\frac{a(\cdot,t)}{\widetilde a}-1\right\|_{H^1}
\|\widetilde a\theta_x\|_{H^1}.
\]
Together with \eqref{eq33}, this shows that the map
\[
\frac{a(\cdot,t)}{\widetilde a}-1
\longmapsto
\left[U\longmapsto
\left(0,0,\partial_x\!\left[
\left(\frac{a(\cdot,t)}{\widetilde a}-1\right)
\widetilde a\theta_x
\right]\right)\right]
\]
is bounded linear map from $H^1(0,1)$ into
$\mathcal L(\mathcal D,\mathcal H)$.
By \eqref{eq7}, one gets 
$
\frac a{\widetilde a}-1
\in W^{1,\infty}(0,T;H^1(0,1)).
$
Therefore, we have 
$
\mathcal A\in
W^{1,\infty}(0,T;\mathcal L(\mathcal D,\mathcal H)),
$
and, for almost every $t$,
$
\partial_t\mathcal A(t)U
=
\left(0,0,\partial_x\!\left[
\frac{a_t}{\widetilde a}\,\widetilde a\theta_x
\right]\right).
$
Moreover, by \eqref{eq33}, 
\[
\|\partial_t\mathcal A(t)U\|_{\mathcal H}
\le C\left\|\frac{a_t}{\widetilde a}\right\|_{H^1}
       \|\widetilde a\theta_x\|_{H^1}
\le C_T\|U\|_{\mathcal D}.
\]
In the zero-flux case, 
since $(\widetilde a\theta_x)(0)=0$, it follows from Hardy's inequality
\eqref{eq25} that 
$
\left\|\frac{\widetilde a\theta_x}{x}\right\|
\le2\|(\widetilde a\theta_x)_x\|.
$
Let us write 
\[
(\mathcal A(t)-\mathcal A(0))U
=
\left(0,0,
\left(\frac{a(\cdot,t)}{\widetilde a}-1\right)
(\widetilde a\theta_x)_x
+
x\partial_x\!\left(\frac{a(\cdot,t)}{\widetilde a}\right)
\frac{\widetilde a\theta_x}{x}
\right),
\]
together with \eqref{eq24}, we obtain
\[
\begin{aligned}
\|(\mathcal A(t)-\mathcal A(0))U\|_{\mathcal H}
&\le
\left(
\left\|\frac{a(\cdot,t)}{\widetilde a}-1\right\|_\infty
+2\left\|x\partial_x\!\left(
\frac{a(\cdot,t)}{\widetilde a}\right)\right\|_\infty
\right)
\|\widetilde a\theta_x\|_{H^1}.
\end{aligned}
\]
Thus, the map
\[
\left(
\frac{a(\cdot,t)}{\widetilde a}-1,\,
x\partial_x\!\left(\frac{a(\cdot,t)}{\widetilde a}\right)
\right)
\longmapsto
\left[U\longmapsto
\left(0,0,
\left(\frac{a(\cdot,t)}{\widetilde a}-1\right)
(\widetilde a\theta_x)_x
+
x\partial_x\!\left(\frac{a(\cdot,t)}{\widetilde a}\right)
\frac{\widetilde a\theta_x}{x}
\right)\right]
\]
is bounded and linear from $L^\infty(0,1)^2$ into
$\mathcal L(\mathcal D,\mathcal H)$. 
By \eqref{eq8}, we have 
$
\left(
\frac a{\widetilde a}-1,\,
x\partial_x\!\left(\frac a{\widetilde a}\right)
\right)
\in W^{1,\infty}(0,T;L^\infty(0,1)^2).
$
Thus, 
$
\mathcal A\in
W^{1,\infty}(0,T;\mathcal L(\mathcal D,\mathcal H)).
$
Since $
\partial_t\!\left[
x\partial_x\!\left(\frac a{\widetilde a}\right)\right]
=
x\partial_x\!\left(\frac{a_t}{\widetilde a}\right)
$
in the distributional sense,  for almost every $t$,
$
\partial_t\mathcal A(t)U
=
\left(0,0,
\frac{a_t}{\widetilde a}(\widetilde a\theta_x)_x
+
\left[x\partial_x\!\left(\frac{a_t}{\widetilde a}\right)\right]
\frac{\widetilde a\theta_x}{x}
\right),
$
with
\[
\begin{aligned}
\|\partial_t\mathcal A(t)U\|_{\mathcal H}
\le
\left(
\left\|\frac{a_t}{\widetilde a}\right\|_\infty
+2\left\|x\partial_x\!\left(
\frac{a_t}{\widetilde a}\right)\right\|_\infty
\right)\|\widetilde a\theta_x\|_{H^1}\le C_T\|U\|_{\mathcal D}.
\end{aligned}
\]
Thus, under either boundary condition, the Bochner derivative is
\begin{equation*}
\partial_t\mathcal A(t)U
=
\left(0,0,
\frac{a_t}{\widetilde a}(\widetilde a\theta_x)_x
+\partial_x\!\left(\frac{a_t}{\widetilde a}\right)
 \widetilde a\theta_x
\right)
\end{equation*}
for almost every $t\in(0,T)$.
In particular,
\[
\partial_t\mathcal A
\in L^\infty(0,T;\mathcal L(\mathcal D,\mathcal H))
\]
in the Bochner sense, verifying \textup{(H$_4$)}
of Theorem~\ref{thm3.4}.

By Proposition~\ref{prop3.5}, for each fixed $t\in[0,T]$,
$\mathcal A(t)$ generates a semigroup $S_t(r)$ satisfying
$
\|S_t(r)\|_{\mathcal L(\mathcal H)}\le1,
$ for $r\ge0.$
All these estimates use the same norm on $\mathcal H$. Hence
\[
\begin{aligned}
\|S_{t_k}(s_k)\cdots S_{t_1}(s_1)U\|_{\mathcal H}
&\le\prod_{j=1}^k\|S_{t_j}(s_j)\|_{\mathcal L(\mathcal H)}
\|U\|_{\mathcal H}
\le\|U\|_{\mathcal H}
\end{aligned}
\]
for every $U\in\mathcal H$, $k\ge1$,
$0\le t_1\le\cdots\le t_k\le T$ and $s_j\ge0$.
Thus \textup{(H$_3$)} of Theorem~\ref{thm3.4} holds
with $C=1$ and $p=0$.

\end{proof}

\begin{proof}[Proof of Theorem~\ref{thm2.6}]
Fix $0\le s<T$. Propositions~\ref{prop3.5} and \ref{prop3.6}
verify the hypotheses of Theorem~\ref{thm3.4} for
$\mathcal A(s+\tau)$ on $[0,T-s]$; the initial graph norms are
equivalent by \eqref{eq28}. Thus each
$U_s\in\mathcal D$ gives a unique strong solution. Dissipativity implies directly  
$E'=-D$, and hence
\eqref{eq14}. By linearity, two strong solutions satisfy
\[
\frac12\|U(t)-\widehat U(t)\|_{\mathcal H}^2
+\int_s^t\!\int_0^1a|(\theta-\widehat\theta)_x|^2\dd x\dd\tau
=\frac12\|U(s)-\widehat U(s)\|_{\mathcal H}^2.
\]

For $U_s\in\mathcal H$, choose $U_s^n\in\mathcal D$ such that
$U_s^n\to U_s$ in $\mathcal H$ as $n\to \infty$. The corresponding regular solutions
satisfy
\[
\begin{aligned}
\sup_{s\le t\le T}\|U^n(t)-U^\ell(t)\|_{\mathcal H}
\le\|U_s^n-U_s^\ell\|_{\mathcal H},\quad
\int_s^T\!\int_0^1a|(\theta^n-\theta^\ell)_x|^2\dd x\dd t
\le\frac12\|U_s^n-U_s^\ell\|_{\mathcal H}^2.
\end{aligned}
\]
Thus $U^n$ converges to some $U$ in $C([s,T];\mathcal H)$ as $n\to \infty$. 
Since $\widetilde a\le a/m$,
\[
\begin{aligned}
\int_s^T\|\theta^n-\theta^\ell\|_V^2\dd t
&=\int_s^T\|\theta^n-\theta^\ell\|^2\dd t
+\int_s^T\!\int_0^1\widetilde a
|(\theta^n-\theta^\ell)_x|^2\dd x\dd t\\
&\le\left(T-s+\frac1{2m}\right)
\|U_s^n-U_s^\ell\|_{\mathcal H}^2
\to0\quad\text{as }n,\ell\to\infty.
\end{aligned}
\]
Completeness of $L^2(s,T;V)$ and uniqueness of the limit in
$L^2(s,T;L^2(0,1))$ identify its limit with the thermal component
$\theta$ of $U$. Consequently, as $n\to \infty$,
\[
\theta^n\to\theta
\quad\text{in }L^2(s,T;V),\quad
\int_s^T\!\int_0^1a|\theta_x^n-\theta_x|^2\dd x\dd t
\le M\|\theta^n-\theta\|_{L^2(s,T;V)}^2
\to 0.
\]
Writing $D_n(t)=\int_0^1a|\theta_x^n|^2\dd x$ and using Cauchy-Schwarz, we obtain, as $n\to \infty$,
\begin{equation}\label{eq34}
\begin{aligned}
\|D_n-D\|_{L^1(s,T)}
&\le\int_s^T\!\int_0^1
a\,|\theta_x^n-\theta_x|
\bigl(|\theta_x^n|+|\theta_x|\bigr)\,dx\,dt\\
&\le
\|\sqrt a(\theta_x^n-\theta_x)\|_{L^2}
\bigl(\|\sqrt a\,\theta_x^n\|_{L^2}
      +\|\sqrt a\,\theta_x\|_{L^2}\bigr)
\to0,
\end{aligned}
\end{equation}
where the $L^2$ norms are over $(s,T)\times(0,1)$.
Convergence in $C([s,T];\mathcal H)$ also illustrates uniform convergence
of $E_n$ to $E$. Together with \eqref{eq34}, this passes
the regular energy identity to the limit for every $t\in[s,T]$.
Thus $E$ is absolutely continuous and $E'=-D$ almost everywhere.
The difference estimate proves independence of the approximation,
uniqueness in the class of Definition~\ref{def2.5}, and
\eqref{eq15}.

Since  $U^n\to U$ in $C([s,T];\mathcal H)$ and $
a\theta_x^n\to a\theta_x$  in $L^2(s,T;L^2(0,1))$ as $n\to \infty$, this allows
passage to the limit in the distributional equations, giving
\[
v_t=(u_x-\kappa\theta)_x,\quad
\theta_t=(a\theta_x-\kappa v)_x.
\]
Passing to the limit in the identity for $u^n$ also shows
\[
u(t)=u_s+\int_s^tv(\tau)\dd\tau
\quad\text{in }L^2(0,1).
\]
Thus $u_t=v$ and the original evolution equations hold in distributions.

In the Dirichlet case, since $u\in C([s,T];H_0^1(0,1))$,
\[
u(0,t)=u(1,t)=0\quad t\in[s,T].
\]
Moreover, $\theta\in L^2(s,T;V)$ shows
\[
\theta(0,t)=\theta(1,t)=0
\quad\text{for almost every }t\in(s,T).
\]

In the zero-flux case, the same argument gives the mechanical Dirichlet conditions
for every time and $\theta(1,t)=0$ for almost every time.
For the degenerate endpoint, integrating the thermal equation together with $u_t=v$ yields
\begin{equation*}
\partial_x\left(
 \int_s^ta(\cdot,\tau)\theta_x(\cdot,\tau)\dd\tau
\right)
=\theta(t)-\theta(s)+\kappa\bigl(u_x(t)-u_x(s)\bigr)
\end{equation*}
in spatial distributions. This identity is obtained by integrating the regularized system and passing to the limit; continuity in time as spatial distributions ensures validity for every $t\in[s,T]$. Since $a\theta_x\in L^2(s,T;L^2(0,1))$, we have
$\int_s^ta(\cdot,\tau)\theta_x(\cdot,\tau)\dd\tau\in C([s,T];L^2(0,1))$.
Thanks to the identity above and the continuity of $\theta,u_x$ in $L^2(0,1)$,
$
\partial_x\left(
 \int_s^ta(\cdot,\tau)\theta_x(\cdot,\tau)\dd\tau
\right)\in C([s,T];L^2(0,1)).
$
Hence
\[
\int_s^{\,\cdot}a\theta_x\dd\tau
\in C([s,T];H^1(0,1)).
\]
Furthermore, it follows from Cauchy--Schwarz in time and $a^2/\widetilde a\le Ma$
that 
\[
\begin{aligned}
\int_0^1
\frac{\left|\int_s^ta\theta_x\dd\tau\right|^2}
     {\widetilde a(x)}\dd x
&\le(t-s)\int_s^t\!\int_0^1
 \frac{a^2}{\widetilde a}|\theta_x|^2\dd x\dd\tau
 \le M(t-s)\int_s^tD(\tau)\dd\tau<\infty.
\end{aligned}
\]
Because $\int_0^11/\widetilde a\dd x=\infty$, the continuous
representative of this integrated flux must vanish at zero,
by the argument of Lemma~\ref{lem3.2}. Therefore
\begin{equation*}
\left[
\int_s^ta(\cdot,\tau)\theta_x(\cdot,\tau)\dd\tau
\right]_{x=0}=0,
\quad t\in[s,T].
\end{equation*}
This is the natural zero-flux condition for energy data, expressed through the trace
of the integrated flux.

The construction applies on every finite interval, and uniqueness
identifies the solutions on overlapping intervals. Thus the solution
extends to $[s,\infty)$.  Defining
$\mathcal U(t,s)U_s=U(t)$, we obtain
\[
\mathcal U(s,s)=I,\quad
\mathcal U(t,r)\mathcal U(r,s)=\mathcal U(t,s),\quad
\|\mathcal U(t,s)\|_{\mathcal L(\mathcal H)}\le1,
\quad s\le r\le t.
\]
Linearity and strong continuity follow from the regular evolution
family and its continuous extension to $\mathcal H$.

\end{proof}

\section{Uniform energy decay}\label{sec:decay}

All auxiliary spatial weights and averaging kernels in this section
are fixed independently of the conductivity. The final decay
constants depend only on $\kappa,M_1,M_2,\widetilde M$
and these fixed choices. Local graph-norm bounds are used to
justify calculations for each strong solution, but do not enter
the final decay constants.

\begin{lemma}
\label{lem4.1}
Assume Hypotheses~\ref{hyp2.2}-\ref{hybound}, and let
$U=(u,v,\theta)$ be a strong solution on a compact time interval $I$. 
Then the functions
\[
t\longmapsto\int_0^1vP\dd x,\quad
t\longmapsto\int_0^1uv\dd x,\quad
t\longmapsto\int_0^1xu_xv\dd x
\]
belong to $C^1(I)$, where $P$ is obtained from $\theta$
by a fixed bounded linear operator on $L^2(0,1)$.
Moreover,
\begin{equation}\label{eq35}
\lim_{x\downarrow0}\sup_{t\in I}
\left\{
x|\theta(x,t)|^2+x|u_x(x,t)|^2+x|v(x,t)|^2
+|(a\theta_x)(x,t)\theta(x,t)|
\right\}=0.
\end{equation}
\end{lemma}

\begin{proof}
Since $U\in C^1(I;\mathcal H)\cap C(I;\mathcal D)$ and $u_t=v$, we obtain 
$
u\in C^1(I;H_0^1(0,1)),~ v,\theta\in C^1(I;L^2(0,1)).
$
Using $\widetilde a=a(\cdot,0)$ and dissipativity, one has 
\[
\int_0^1\widetilde a\,\theta_x^2\dd x
=-\langle\mathcal A(0)U,U\rangle_{\mathcal H}
\le\|\mathcal A(0)U\|_{\mathcal H}\|U\|_{\mathcal H}.
\]
Together with  \eqref{eq30}, this yields
\[
\|v\|_{H_0^1}+\|\theta\|_V
+\|u_x-\kappa\theta\|_{H^1}
+\|\widetilde a\theta_x\|_{H^1}
\le C_T\|U\|_{\mathcal D}.
\]
The corresponding linear maps are therefore continuous, thus \(\theta\in C(I;V)\), $
v\in C(I;H_0^1(0,1))$, $u_x-\kappa\theta,\ \widetilde a\theta_x\in C(I;H^1(0,1))$.  Using $a\theta_x=(a/\widetilde a)(\widetilde a\theta_x)$, and using the $H^1$ product estimate in the Dirichlet case and Lemma \ref{lem3.3} in the zero-flux case, we obtain
\[
\begin{aligned}
\|a(t)\theta_x(t)-a(s)\theta_x(s)\|_{H^1}
\le
C_T\|\widetilde a\theta_x(t)-\widetilde a\theta_x(s)\|_{H^1}+\left\|
\frac{a(t)-a(s)}{\widetilde a}\,
\widetilde a\theta_x(s)
\right\|_{H^1}
\to0
\end{aligned}
\]
as $t\to s$, by $\widetilde a\theta_x\in C(I;H^1)$ and the
coefficient regularity in Hypothesis~\ref{hyp2.3}. 
Thus $a\theta_x\in C(I;H^1(0,1))$. 
For $0<\varepsilon<1$, positivity of $\widetilde a$ on
$[\varepsilon,1]$ gives
$
\|\theta_x\|_{L^2(\varepsilon,1)}
\le
\frac{\|\sqrt{\widetilde a}\,\theta_x\|}
     {\sqrt{\min_{[\varepsilon,1]}\widetilde a}},
$
hence $\theta\in C(I;H^1(\varepsilon,1))$.
The identity
$u_x=(u_x-\kappa\theta)+\kappa\theta$ then gives
$u\in C(I;H^2(\varepsilon,1))$. Since $u,v,P,u_x\in C^1(I;L^2(0,1))$, with $u_t=v$
and $(u_x)_t=v_x$, the Hilbert-space product rule shows that
the three integrals belong to $C^1(I)$. The same argument
applies when their integrands are multiplied by fixed smooth
spatial weights. Moreover, for every $\varepsilon\in(0,1)$,
the local regularity
$u\in H^2(\varepsilon,1)$ and $v,\theta\in H^1(\varepsilon,1)$
justifies spatial integration by parts on $(\varepsilon,1)$,
where all resulting boundary traces are well defined.

Since $\theta\in C(I;V)\subset C(I;L^2),$ and $I$ is compact,
the set $\{\theta(t):t\in I\}$ is compact in $L^2(0,1)$ and $
\sup_{t\in I}\|\theta(t)\|_V<\infty.$ 
A finite-covering argument and absolute continuity of the integral
therefore imply 
$\sup_{t\in I}\|\theta(t)\|_{L^2(0,x)}\to0$
as $x\downarrow0$.
Moreover, boundedness in $V$ and the weight assumption yield
\[
\sup_{t\in I}\|x\theta_x(t)\|^2
\le
\left(\sup_{0<x\le1}\frac{x^2}{\widetilde a(x)}\right)
\sup_{t\in I}\|\theta(t)\|_V^2<\infty.
\]
For fixed $t\in I$, integrating
$\partial_y(y|\theta(y,t)|^2)
=|\theta(y,t)|^2+2y\theta(y,t)\theta_x(y,t)$
over $(\varepsilon,x)$ and letting $\varepsilon\downarrow0$,
using \eqref{eq21} with $z=\theta$, yields
\[
x|\theta(x,t)|^2
=\int_0^x
\bigl(|\theta(y,t)|^2+2y\theta(y,t)\theta_x(y,t)\bigr)\dd y
\le\|\theta(t)\|_{L^2(0,x)}^2
+2\|\theta(t)\|_{L^2(0,x)}\|x\theta_x(t)\|,\]
which tends to zero uniformly as $x\downarrow0$. By the boundedness of $u_x-\kappa\theta$ and $v$ in $H^1(0,1)$, we obtain 
\[
\begin{aligned}
x|u_x(x,t)|^2
\le2x\|u_x(t)-\kappa\theta(t)\|_\infty^2
+2\kappa^2x|\theta(x,t)|^2,\quad 
x|v(x,t)|^2&\le x\|v(t)\|_\infty^2.
\end{aligned}
\]
Thus both mechanical boundary terms also vanish uniformly as $x\downarrow0$.

In the Dirichlet case,
using $\theta(0,t)=0$, $1/\widetilde a\in L^1(0,1)$ and Cauchy-Schwarz, then 
\[
|\theta(x,t)|^2=\left|\int_0^x\theta_x(y,t)\,dy\right|^2
\le
\left(\int_0^x\frac{\dd y}{\widetilde a(y)}\right)
\int_0^1\widetilde a(y)|\theta_x(y,t)|^2\dd y.
\]
Hence $\sup_{t\in I}|\theta(x,t)|\to0$ as $x\downarrow0$.
Since $a\theta_x$ is bounded in $H^1(0,1)$, it follows that
\[
\sup_{t\in I}|(a\theta_x)(x,t)\theta(x,t)|
\le
\sup_{t\in I}\|a\theta_x(t)\|_\infty
\sup_{t\in I}|\theta(x,t)|
\to0\quad\text{as } x\downarrow0.
\]
In the zero-flux case, it follows from $(a\theta_x)(0,t)=0$ and
Cauchy--Schwarz that 
$
|(a\theta_x)(x,t)|^2
=\left|\int_0^x(a\theta_x)_x(y,t)\dd y\right|^2
\le x\|(a\theta_x)_x(t)\|^2.
$
Consequently,
\[
\sup_{t\in I}|(a\theta_x)(x,t)\theta(x,t)|^2
\le
\left(\sup_{t\in I}\|(a\theta_x)_x(t)\|^2\right)
\left(\sup_{t\in I}x|\theta(x,t)|^2\right)
\to 0
\quad\text{as }x\downarrow0.
\]
 Together with the preceding estimates, 
this proves \eqref{eq35} in both boundary cases.

\end{proof}

Now, we recall the following well-known result \cite[Theorem 8.1]{Komornik1994}.
\begin{lemma}\label{lem4.2}
Assume that $E : [0,+\infty) \to [0,+\infty)$ is a nonincreasing function and that there is a constant $\mathcal{M}>0$ such that
\begin{equation*}
\int_t^\infty E(s)\,ds \leq \mathcal{M} E(t), \quad \forall t\in[0,+\infty).
\end{equation*}
Then we have
\begin{equation*}
E(t) \leq E(0) e^{1-t/\mathcal{M}}, \quad \forall t\in[\mathcal{M},+\infty).
\end{equation*}
\end{lemma}

\subsection{Dirichlet boundary condition}
\label{sec:dirichlet}

The argument combines boundary compensation with time averaging:
we first eliminate the endpoint heat flux, then estimate the
resulting low-order mechanical terms on a fixed time interval,
and finally absorb them to obtain an energy contraction. 

\begin{proposition}[Interior-average compensation]
\label{prop4.3}
Assume Hypotheses~\ref{hyp2.2}-\ref{hybound}, in the Dirichlet case. 
Fix $\kappa>0$ and a nonnegative function
$\rho\in C_c^\infty(0,1)$ satisfying
$\int_0^1\rho\dd x=1$.
There exists an instantaneous quadratic functional
${\mathcal L_1}$, equivalent to $E$, such that every strong
solution satisfies
\begin{equation}\label{eq36}
{\mathcal L_1}'(t)
\le-cE(t)+C_{\rm dec}\bigl(\|u(t)\|^2+|m_\rho(t)|^2\bigr),
\quad
m_\rho(t)=\int_0^1\rho(x)v(x,t)\dd x.
\end{equation}
The equivalence constants and $c,C_{\rm dec}>0$ depend only on
$\kappa,M_1,M_2,\widetilde{M}$ and the fixed spatial weights, independently of the time-regularity bounds.
\end{proposition}

\begin{proof}
\emph{Step 1: removal of the endpoint flux.}
Define
\begin{equation*}
P(x,t)=\int_0^x\theta(y,t)\dd y
-\int_0^1\rho(z)\int_0^z\theta(y,t)\dd y\dd z.
\end{equation*}
Then $P_x=\theta$, $\int_0^1\rho P\dd x=0$. 
By Cauchy--Schwarz, $\rho\ge0$, and $\int_0^1\rho\dd x=1$,
\[
\begin{aligned}
|P(x,t)|
&\le
\left|\int_0^x\theta(y,t)\dd y\right|
+\int_0^1\rho(z)
\left|\int_0^z\theta(y,t)\dd y\right|\dd z\\
&\le
\left(\sqrt{x}+\int_0^1\rho(z)\sqrt{z}\dd z\right)
\|\theta(t)\|
\le2\|\theta(t)\|.
\end{aligned}
\]
Using $P_x=\theta$, we conclude that
\[
\|P(t)\|_{H^1}^2
=\|P(t)\|^2+\|\theta(t)\|^2
\le5\|\theta(t)\|^2.
\]
The thermal equation and $P_x=\theta$ show 
$\partial_x(P_t-a\theta_x+\kappa v)=0$.
Taking the $\rho$-average and using
$\int_0^1\rho P_t\dd x=0$ and $\int_0^1\rho\dd x=1$,
we obtain
\begin{equation}\label{eq37}
P_t=a\theta_x-\int_0^1\rho a\theta_x\dd x
-\kappa(v-m_\rho).
\end{equation}
Thus the endpoint flux is replaced by an \emph{interior average} $\int_0^1\rho a\theta_x\dd x$,
which is controlled by dissipation:
\[
\left|\int_0^1\rho a\theta_x\dd x\right|
\le
\left(\int_0^1\rho^2\dd x\right)^{1/2}
\left(\int_0^1a^2|\theta_x|^2\dd x\right)^{1/2}
\le\|\rho\|\sqrt{\widetilde{M}D}.
\]
Here we use Cauchy–Schwarz and \eqref{eq20} with $z=\theta$.

\emph{Step 2: the interior velocity estimate.} 
Choose $\chi\in C_c^\infty(0,1)$ with $0\le\chi\le1$ and
$\chi=1$ on $\operatorname{supp}\rho$.  
Define \[F=\int_0^1\chi vP\dd x.\]
 Thus, by Lemma~\ref{lem4.1},
\eqref{eq37}  and integration by parts, 
\begin{equation*}
\begin{aligned}
F'={}&\int_0^1\chi v_tP\,dx+\int_0^1\chi vP_t\,dx\\
={}&-\kappa\int_0^1\chi v^2\dd x
-\int_0^1\chi'(u_x-\kappa\theta)P\dd x
-\int_0^1\chi(u_x-\kappa\theta)\theta\dd x\\
&+\int_0^1\chi va\theta_x\dd x
-\left(\int_0^1\rho a\theta_x\dd x\right)
 \left(\int_0^1\chi v\dd x\right)
+\kappa m_\rho\int_0^1\chi v\dd x.
\end{aligned}
\end{equation*}
Using $\|P\|\le2\|\theta\|$, \eqref{eq20} with $z=\theta$ and Cauchy--Schwarz arrives at 
\[
\left|F'+\kappa\int_0^1\chi v^2\dd x\right|
\le C\bigl((\|u_x\|+\|v\|)\sqrt D+D
+|m_\rho|\|v\|\bigr).
\]
Young's inequality therefore yields, for every $\varepsilon>0$,
\begin{equation}\label{eq38}
F'\le-\kappa\int_0^1\chi v^2\dd x
+\varepsilon E+C_\varepsilon(D+|m_\rho|^2),
\end{equation}
where $C_\varepsilon$ depends only on
$\varepsilon,\kappa,M_2,\widetilde{M}$ and the fixed weights.
Moreover,
\begin{equation}
\label{eq39}
|F|\le\|v\|\|P\|
\le2\|v\|\|\theta\|
\le\|v\|^2+\|\theta\|^2\le2E.
\end{equation}

\emph{Step 3: the mechanical multiplier.} 
Define $h(x)=x-\int_0^x\rho(y)\dd y$, then 
\begin{equation}
\label{eq40}
h(0)=h(1)=0, \quad h'=1-\rho, \quad 
|h(x)|\le(1+\|\rho\|_\infty)x.
\end{equation}
Further, using the first in \eqref{eq20} with $z=\theta$, one has 
\begin{equation}
\label{eq41}
\|h\theta_x\|^2\le(1+\|\rho\|_\infty)^2M_1D.
\end{equation}

Define 
\begin{equation*}
G(t)=\int_0^1v\left(hu_x+\frac12\rho u\right)\,dx
\end{equation*}
Using $u_t=v$ and $v_t=(u_x-\kappa\theta)_x$, we obtain
\[
G'
=\int_0^1h(v_xv+u_xv_t)\dd x
+\frac12\int_0^1\rho(v^2+uv_t)\dd x.
\]
For the first term, perform spatial integration by parts
on $(x,1)$ and then integrate in time. The boundary terms
vanish uniformly as $x\downarrow0$ by \eqref{eq40} and
\eqref{eq35}, and the mixed volume term is integrable by
\eqref{eq41}. Passing to the limit and differentiating the
resulting time-integrated identity yields
\[
\int_0^1h(v_xv+u_xv_t)\dd x
=-\frac12\int_0^1(1-\rho)(v^2+u_x^2)\dd x
-\kappa\int_0^1h\theta_xu_x\dd x.
\]
Since $\rho$ has compact support, integration by parts and \(uu_x=\frac12(u^2)_x\) results in 
\[
\begin{aligned}
\frac12\int_0^1\rho(v^2+uv_t)\dd x
={}&\frac12\int_0^1\rho v^2\dd x
-\frac12\int_0^1\rho u_x^2\dd x
+\frac14\int_0^1\rho''u^2\dd x\\
&+\frac\kappa2\int_0^1\rho u_x\theta\dd x
+\frac\kappa2\int_0^1\rho'u\theta\dd x.
\end{aligned}
\]
Combining the above two identities, we obtain
\[
\begin{aligned}
G'&=-\frac12(\|v\|^2+\|u_x\|^2)
+\int_0^1\rho v^2\dd x
+\frac14\int_0^1\rho''u^2\dd x\\
&\quad-\kappa\int_0^1h\theta_xu_x\dd x
+\frac\kappa2\int_0^1\rho u_x\theta\dd x
+\frac\kappa2\int_0^1\rho'u\theta\dd x.
\end{aligned}
\]
By \eqref{eq20} with $z=\theta$, \eqref{eq41}, $\|u\|\le\|u_x\|$
and Young's inequality,
\[
\begin{aligned}
&\left|\kappa\int_0^1h\theta_xu_x\dd x\right|
+\left|\frac\kappa2\int_0^1\rho u_x\theta\dd x\right|
+\left|\frac\kappa2\int_0^1\rho'u\theta\dd x\right|\\
&\quad\le
\kappa(1+\|\rho\|_\infty)\sqrt{M_1D}\,\|u_x\|
+\frac\kappa2\|\rho\|_\infty\sqrt{M_2 D}\,\|u_x\|
+\frac\kappa2\|\rho'\|_\infty\sqrt{M_2 D}\,\|u_x\|\\
&\quad\le\frac3{16}\|u_x\|^2+CD.
\end{aligned}
\]
Since $\rho\le\|\rho\|_\infty\chi$, it follows that
\[
G'\le-\frac12\|v\|^2-\frac5{16}\|u_x\|^2
+C_0\int_0^1\chi v^2\dd x
+\widetilde C_0D+C_2\|u\|^2,
\]
where $C_0=\|\rho\|_\infty$,
$C_2=\|\rho''\|_\infty/4$, and
$\widetilde C_0$ depends only on
$\kappa,M_1,M_2$ and $\rho$.
Using \eqref{eq20} with $z=\theta$ and setting
$C_1=\widetilde C_0+M_2/4$, we obtain
$
-\frac12\|v\|^2-\frac5{16}\|u_x\|^2+\widetilde C_0D
\le-\frac12E+C_1D.
$
Thus, we get 
\begin{equation}\label{eq42}
G'\le-\frac12E+C_0\int_0^1\chi v^2\dd x
+C_1D+C_2\|u\|^2.
\end{equation}
Moreover,
\begin{equation}
\label{eq43}
|G|
\le\|h\|_\infty\|u_x\|\|v\|
+\frac12\|\rho\|_\infty\|u\|\|v\|
\le C_3E,
\end{equation}
with $
C_3=\|h\|_\infty+\frac12\|\rho\|_\infty.$
\medskip

\emph{Step 4: combination of the estimates.} 
Define 
\[\mathcal L_1=N_1E+G+N_2F.\] 
Combining $E'=-D$ with
\eqref{eq38} and
\eqref{eq42}, we obtain
\[
\begin{aligned}
\mathcal L_1'
\le
\left(-\frac12+N_2\varepsilon\right)E
+(C_0-N_2\kappa)\int_0^1\chi v^2\dd x+(-N_1+C_1+N_2C_\varepsilon)D
+C_2\|u\|^2+N_2C_\varepsilon|m_\rho|^2.
\end{aligned}
\]
Choose successively
$
N_2=\frac{C_0+1}{\kappa},~
\varepsilon=\frac1{4N_2},~
N_1\ge\max\{2(C_3+2N_2),\,C_1+N_2C_\varepsilon+1\}.
$
Since $|G|+N_2|F|\le(C_3+2N_2)E\le N_1E/2$ by \eqref{eq39} and \eqref{eq43}, we have
\begin{equation}\label{eq44}
\frac {N_1}{2}E\le\mathcal L_1\le\frac{3N_1}{2}E.
\end{equation}
With these choices, further, 
\[
\begin{aligned}
\mathcal L_1'
&\le-\frac14E-\int_0^1\chi v^2\dd x-D
+C_2\|u\|^2+N_2C_\varepsilon|m_\rho|^2\le-\frac14E+C_{\rm dec}
\bigl(\|u\|^2+|m_\rho|^2\bigr),
\end{aligned}
\]
where $C_{\rm dec}=C_2+N_2C_\varepsilon$.
This proves \eqref{eq36}.

\end{proof}

The remainder $\|u(t)\|^2+|m_\rho(t)|^2$ in \eqref{eq36} is controlled by
time averaging: the thermal equation estimates the averaged velocity
in a negative norm, spatial interpolation yields an $L^2$ estimate,
and the mechanical equation recovers the displacement. Averaging
the actual thermal flux avoids differentiating $a$ in time.

\begin{lemma}[Quantitative recovery of the low-order mechanical terms]
\label{lem4.4}
Assume Hypotheses~\ref{hyp2.2}-\ref{hybound}, in the Dirichlet case,
with $\kappa>0$, and let $\rho$ be as in 
Proposition~\ref{prop4.3}.
For every $\epsilon>0$, there exists $C_\epsilon>0$ such that
every energy solution satisfies, for all $s\ge0$ and $T\ge1$,
\begin{equation}\label{eq45}
\begin{aligned}
\int_s^{s+T}\bigl(\|u(t)\|^2+|m_\rho(t)|^2\bigr)\dd t
\le
\epsilon\int_s^{s+T}E(t)\dd t+C_\epsilon\int_s^{s+T}D(t)\dd t
+C_\epsilon E(s).
\end{aligned}
\end{equation}
The constant depends only on
$\epsilon,\kappa,M_2,\widetilde{M}$, the fixed weight $\rho$,
and the averaging kernel. 
\end{lemma}

\begin{proof}
Let $U$ be an energy solution, fix $s\ge0$ and $T\ge1$,
and set $I=[s,s+T]$.
Theorem~\ref{thm2.6} gives
$
u\in C(I;H_0^1(0,1))\cap C^1(I;L^2(0,1)),$ $v,\theta\in C(I;L^2(0,1)),
$ with $u_t=v$, $a\theta_x\in L^2(I;L^2(0,1))$, and the
mechanical and thermal equations holding in distributions.
We use only these properties, the energy identity and the
spatial bounds in \eqref{eq20}.

\emph{Step 1: a spatial interpolation estimate.} 
Equip $H_0^1(0,1)$ with the norm $\|\varphi_x\|$ and
$H^{-1}(0,1)$ with the corresponding dual norm.
For $w\in H_0^1(0,1)$, set $\overline w=\int_0^1w\dd x$. Since $\int_0^1\varphi_x\dd x=0$ for $\varphi\in H_0^1(0,1)$,
and every zero-mean $g\in L^2(0,1)$ is the derivative of
$\varphi(x)=\int_0^xg(y)\dd y\in H_0^1(0,1)$, we have
\[
\begin{aligned}
\|\partial_xw\|_{H^{-1}}
=\sup_{\substack{\varphi\in H_0^1\\\|\varphi_x\|=1}}
\left|\int_0^1(w-\overline w)\varphi_x\dd x\right|
=\sup_{\substack{g\in L^2,\ \int_0^1g\dd x=0\\\|g\|=1}}
\left|\int_0^1(w-\overline w)g\dd x\right|
=\|w-\overline w\|.
\end{aligned}
\]
Moreover, $w(0)=0$ shows $|w(x)|^2\le x\|w_x\|^2$.
Hence, for $0<\delta<1$,
\[
\delta|\overline w|^2
\le2\int_0^\delta|w|^2\dd x
+2\|w-\overline w\|^2
\le\delta^2\|w_x\|^2+2\|w-\overline w\|^2.
\]
Combining this with
$\|w\|^2=\|w-\overline w\|^2+|\overline w|^2$ yields
\begin{equation}\label{eq46}
\|w\|^2
\le\delta\|w_x\|^2
+(1+2/\delta)\|\partial_xw\|_{H^{-1}}^2.
\end{equation}

\emph{Step 2: temporal averaging and control of the averaged velocity.}
Choose a nonnegative $\eta\in C_c^\infty(-1,1)$ with
$\int_{\mathbb R}\eta\dd r=1$. For $0<h<1/4$, define
\[
\eta_h(r)=\frac1h\eta\left(\frac rh\right),\quad
I_h=[s+h,s+T-h],\quad
I_{2h}=[s+2h,s+T-2h].
\]
Extend $u,v,\theta$ and $a\theta_x$ by zero outside $I$ and
use convolution in time to define
\[
\begin{gathered}
u_h=\eta_h*u:=\int_{-h}^{h}\eta_h(r)u(t-r)\,dr,\quad v_h=\eta_h*v,\quad
\theta_h=\eta_h*\theta,\quad
(a\theta_x)_h=\eta_h*(a\theta_x),\\
u_{hh}=\eta_h*\eta_h*u,\quad
\theta_{hh}=\eta_h*\eta_h*\theta.
\end{gathered}
\]
Writing $c_\eta=\|\eta'\|_{L^1(\mathbb R)}$, we have
$
\|\eta_h\|_{L^1(\mathbb R)}=1,~\|\eta_h'\|_{L^1(\mathbb R)}=\frac{c_\eta}{h}.
$
Averaged equations are used on $I_h$ for a single convolution
and on $I_{2h}$ for a double convolution, so the kernels avoid
the temporal endpoints and commute with time differentiation.

Since $u\in C(I;H_0^1(0,1))$, convolution in time makes
$u_h$ smooth with values in $H_0^1(0,1)$ on $I_h$.
The distributional identity $u_t=v$ therefore yields
\[
v_h=\partial_tu_h\in H_0^1(0,1),
\quad (v_h)_x=\eta_h'*u_x
\quad\text{on }I_h.
\]
Young's convolution inequality and $\|u_x\|^2\le2E$ yield
\begin{equation}\label{eq47}
\int_{I_h}\|(v_h)_x\|^2\dd t
\le\|\eta_h'\|_{L^1}^2\int_I\|u_x\|^2\dd t
=\frac{c_\eta^2}{h^{2}}\int_I\|u_x\|^2\dd t
\le\frac{2c_\eta^2}{h^{2}}\int_I E(t)\dd t.
\end{equation}
Convolving the distributional thermal equation in time gives
\[
\kappa(v_h)_x
=\partial_x((a\theta_x)_h)-\eta_h'*\theta
\quad\text{in }H^{-1}(0,1).
\]
Using the dual norm on $H^{-1}(0,1)$ and Poincaré's inequality,
followed by Young's convolution inequality, we obtain
\[
\begin{aligned}
\int_{I_h}\|(v_h)_x\|_{H^{-1}}^2\dd t
&\le\frac2{\kappa^2}\left(
\int_{I_h}\|(a\theta_x)_h\|^2\dd t
+\int_{I_h}\|\eta_h'*\theta\|^2\dd t
\right)\\
&\le\frac2{\kappa^2}\left(
\int_I\|a\theta_x\|^2\dd t
+\frac{c_\eta^2}{h^{2}}\int_I\|\theta\|^2\dd t
\right).
\end{aligned}
\]
Since $\|a\theta_x\|^2\le \widetilde{M}D$ and
$\|\theta\|^2\le M_2 D$ by \eqref{eq20}, it follows that
\begin{equation}\label{eq48}
\int_{I_h}\|(v_h)_x\|_{H^{-1}}^2\dd t
\le C_h\int_I D(t)\dd t,
\quad
C_h=\frac2{\kappa^2}
\bigl(\widetilde{M}+M_2 \frac{c_\eta^2}{h^{2}}\bigr).
\end{equation}
Applying \eqref{eq46} to
$v_h(t)\in H_0^1(0,1)$ and combining
\eqref{eq47}--\eqref{eq48}, we find
\[
\int_{I_h}\|v_h\|^2\dd t
\le\frac{2\delta c_\eta^2}{h^{2}}\int_I E(t)\dd t
+(1+2/\delta)C_h\int_I D(t)\dd t.
\]
For any $\tau>0$, choose
$\delta=\min\{1/2,\tau h^2/(2c_\eta^2)\}$. Then
\begin{equation}\label{eq49}
\int_{I_h}\|v_h\|^2\dd t
\le\tau\int_I E(t)\dd t
+C_{\tau,h}\int_I D(t)\dd t,
\quad
C_{\tau,h}=(1+2/\delta)C_h.
\end{equation}

\emph{Step 3: recovery of the velocity average.}
Testing the distributional mechanical equation against
$\rho\in C_c^\infty(0,1)$ gives the following identity for
$m_\rho'$.
Its right-hand side is continuous in time because $u_x$ and
$\theta$ are continuous in $L^2(0,1)$, so $m_\rho\in C^1(I)$:
\[
m_\rho'
=\langle v_t,\rho\rangle_{H^{-1},H_0^1}
=-\int_0^1\rho'(u_x-\kappa\theta)\dd x.
\]
Therefore, by Cauchy-Schwarz, one has 
\[
|m_\rho'|^2\le\|\rho'\|^2\|u_x-\kappa\theta\|^2
\le2\|\rho'\|^2\bigl(\|u_x\|^2+\kappa^2\|\theta\|^2\bigr)
\le C_\rho E, \quad
C_\rho=4\|\rho'\|^2\max\{1,\kappa^2\}.
\]
Since $\eta_h\ge0$ and $\int_{-h}^h\eta_h(r)\dd r=1$, we have
\[
m_\rho(t)-(m_\rho)_h(t)
=\int_{-h}^h\eta_h(r)
\bigl(m_\rho(t)-m_\rho(t-r)\bigr)\dd r
=\int_{-h}^h\eta_h(r)
r\int_0^1m_\rho'(t-\sigma r)\dd\sigma\dd r.
\]
For $|r|\le h$ and $0\le\sigma\le1$, we have
$I_h-\sigma r\subset I$, and hence
\[
\int_{I_h}|m_\rho'(t-\sigma r)|^2\dd t
=\int_{I_h-\sigma r}|m_\rho'(\tau)|^2\dd\tau
\le\int_I|m_\rho'(\tau)|^2\dd\tau.
\]
Therefore, Jensen's inequality and Cauchy--Schwarz yield
\begin{equation}\label{1110}
\begin{aligned}
\int_{I_h}|m_\rho-(m_\rho)_h|^2\dd t
&=\int_{I_h}
\left|\int_{-h}^h\eta_h(r)
\bigl(m_\rho(t)-m_\rho(t-r)\bigr)\dd r\right|^2\dd t\\
&\le\int_{-h}^h\eta_h(r)
\int_{I_h}|m_\rho(t)-m_\rho(t-r)|^2\dd t\dd r\\
&=\int_{-h}^h\eta_h(r)r^2
\int_{I_h}\left|
\int_0^1m_\rho'(t-\sigma r)\dd\sigma
\right|^2\dd t\dd r\\
&\le\int_{-h}^h\eta_h(r)r^2
\int_0^1\int_{I_h}
|m_\rho'(t-\sigma r)|^2\dd t\dd\sigma\dd r\\
&\le\left(\int_{-h}^h\eta_h(r)r^2\dd r\right)
\int_I|m_\rho'(t)|^2\dd t\\
&\le h^2\int_I|m_\rho'(t)|^2\dd t
\le C_\rho h^2\int_I E(t)\dd t.
\end{aligned}
\end{equation}
Moreover, $(m_\rho)_h=\int_0^1\rho v_h\dd x$, so
Cauchy--Schwarz and \eqref{eq49} yield
\[
\int_{I_h}|(m_\rho)_h|^2\dd t
\le\|\rho\|^2\left(
\tau\int_I E(t)\dd t
+C_{\tau,h}\int_I D(t)\dd t
\right).
\]
Combining these estimates with
$|m_\rho|^2\le2|m_\rho-(m_\rho)_h|^2+2|(m_\rho)_h|^2$,
we obtain
\begin{equation}\label{eq50}
\int_{I_h}|m_\rho|^2\dd t
\le(2C_\rho h^2+2\|\rho\|^2\tau)\int_I E(t)\dd t
+2\|\rho\|^2C_{\tau,h}\int_I D(t)\dd t.
\end{equation}

\emph{Step 4: recovery of the displacement.}
The double kernel $\eta_h*\eta_h$ is nonnegative, has unit mass
and is supported in $(-2h,2h)$.
Applying the argument of \eqref{1110} to $u$ with this kernel
and using $u_t=v$, we obtain
\begin{equation}\label{eq51}
\int_{I_{2h}}\|u-u_{hh}\|^2\dd t
\le4h^2\int_I\|v\|^2\dd t
\le8h^2\int_I E(t)\dd t.
\end{equation}
On $I_{2h}$, the averaged mechanical equation reads
\[
(u_{hh})_{xx}=\eta_h*\eta_h*v_t+\kappa(\theta_{hh})_x
=\eta_h'*v_h+\kappa(\theta_{hh})_x
\quad\text{in }H^{-1}(0,1).
\]
Testing against $u_{hh}\in H_0^1(0,1)$ and using
Cauchy--Schwarz and Poincaré's inequality, we obtain
\[
\|(u_{hh})_x\|^2
\le\|\eta_h'*v_h\|\|u_{hh}\|
+\kappa\|\theta_{hh}\|\|(u_{hh})_x\|
\le
\bigl(\|\eta_h'*v_h\|+\kappa\|\theta_{hh}\|\bigr)
\|(u_{hh})_x\|,
\]
and hence
\[
\|u_{hh}\|^2
\le\|(u_{hh})_x\|^2
\le2\|\eta_h'*v_h\|^2+2\kappa^2\|\theta_{hh}\|^2.
\]
For $t\in I_{2h}$, the first convolution uses only values
of $v_h$ on $I_h$. Thus Young's convolution inequality,
\eqref{eq49} and $\|\theta\|^2\le M_2 D$ by \eqref{eq20} yield
\[
\begin{aligned}
\int_{I_{2h}}\|u_{hh}\|^2\dd t
&\le\frac{2c_\eta^2}{h^2}\int_{I_h}\|v_h\|^2\dd t
+2\kappa^2\int_I\|\theta\|^2\dd t\\
&\le\frac{2c_\eta^2\tau}{h^2}\int_I E(t)\dd t
+\left(\frac{2c_\eta^2C_{\tau,h}}{h^2}
+2\kappa^2M_2\right)\int_I D(t)\dd t.
\end{aligned}
\]
Combining this with \eqref{eq51} and
$\|u\|^2\le2\|u-u_{hh}\|^2+2\|u_{hh}\|^2$, we conclude that
\begin{equation}\label{eq52}
\begin{aligned}
\int_{I_{2h}}\|u\|^2\dd t
\le
\left(16h^2+\frac{4c_\eta^2\tau}{h^2}\right)
\int_I E(t)\dd t
+\left(\frac{4c_\eta^2C_{\tau,h}}{h^2}
+4\kappa^2M_2\right)\int_I D(t)\dd t.
\end{aligned}
\end{equation}

\emph{Step 5: parameter choices and temporal boundary strips.}
Combining \eqref{eq50} and \eqref{eq52}
on $I_{2h}\subset I_h$, we obtain
\begin{equation}
\label{11101}
\begin{aligned}
\int_{I_{2h}}(\|u\|^2+|m_\rho|^2)\dd t
\le
C\left[h^2+\tau\left(1+\frac{1}{h^{2}}\right)\right]\int_I E(t)\dd t
+\widetilde C_{\tau,h}\int_I D(t)\dd t,
\end{aligned}
\end{equation}
where $C$ depends only on $\kappa,\rho,\eta$, and
$\widetilde C_{\tau,h}$ additionally depends on
$M_2,\widetilde{M},\tau,h$.
Choose
\[
h=\min\left\{\frac18,\sqrt{\frac{\epsilon}{2C}}\right\},
\quad
\tau=\frac{\epsilon}{2C\left(1+\frac{1}{h^{2}}\right)},
\]
so that the coefficient of $\int_I E(t)\dd t$ is at most
$\epsilon$. By Poincaré's inequality and Cauchy–Schwarz, one has 
\[
\|u(t)\|^2+|m_\rho(t)|^2
\le \|u_x(t)\|^2+\|\rho\|^2\|v(t)\|^2
\le2(1+\|\rho\|^2)E(t)
\le2(1+\|\rho\|^2)E(s),
\]
which together with $|I\setminus I_{2h}|=4h$, we obtain that 
the remaining strips satisfy
\begin{equation}
\label{11102}
\int_{I\setminus I_{2h}}(\|u\|^2+|m_\rho|^2)\dd t
\le8h(1+\|\rho\|^2)E(s).
\end{equation}
Adding \eqref{11101} and \eqref{11102} proves \eqref{eq45} with
$C_\epsilon\ge
\max\{\widetilde C_{\tau,h},8h(1+\|\rho\|^2)\}$,
independently of $s,T$ and the time-regularity norms.
\end{proof}

\begin{remark}[A common compensation argument]
The interior-average argument of Proposition~\ref{prop4.3}
and Lemma~\ref{lem4.4} also applies to the zero-flux case.
Indeed, \eqref{eq37} uses $P_x=\theta$ and
$\int_0^1\rho P\dd x=0$, and the required mechanical endpoint
limits hold in both cases by Lemma~\ref{lem4.1}.
This gives a common route to exponential decay.
We retain the separate zero-flux proof because it provides
an instantaneous strict Lyapunov functional and explicit
constants without time averaging.
\end{remark}

\begin{proposition}[Integral estimate and decay]
\label{prop4.5}
Assume Hypotheses~\ref{hyp2.2}-\ref{hybound}, in the Dirichlet
case, and let $\kappa>0$. 
There exists $\mathcal M\ge1$, depending only on
$\kappa,M_1,M_2,\widetilde{M}$ and the fixed spatial weights
and averaging kernel, such that every energy solution satisfies
\begin{equation}\label{eq53}
\int_s^{s+T}E(t)\dd t\le \mathcal{M}E(s),
\quad s\ge0,\quad T\ge1.
\end{equation}
Consequently,
\begin{equation}\label{eq54}
E(t)\le
\exp\!\left(1-\frac{t-s}{\mathcal{M}}\right)E(s),
\quad t\ge s\ge0.
\end{equation}
\end{proposition}

\begin{proof}
For a strong solution, \eqref{eq44} and
\eqref{eq36} show
\[
0\le\mathcal L_1\le \hat{c}E,\quad
\mathcal L_1'\le-cE
+C_{\rm dec}(\|u\|^2+|m_\rho|^2),
\]
where $c=1/4$ and $\hat{c}=3N_1/2$. 
Fix $I=[s,s+T]$ with $T\ge1$ and choose
$\epsilon=c/(2C_{\rm dec})$ in \eqref{eq45}.
Integrating the differential inequality, using
$\mathcal L_1(s+T)\ge0$, and applying \eqref{eq45}, we obtain
\[
\begin{aligned}
c\int_I E\dd t
\le \hat{c}E(s)
+C_{\rm dec}\epsilon\int_I E\dd t
+C_{\rm dec}C_\epsilon\int_I D\dd t
+C_{\rm dec}C_\epsilon E(s).
\end{aligned}
\]
Since \eqref{eq14} gives $\int_I D\dd t\le E(s)$,
the chosen value of $\epsilon$ yields
\[
\frac c2\int_I E\dd t
\le(\hat{c}+2C_{\rm dec}C_\epsilon)E(s).
\]
Thus \eqref{eq53} holds with
\[
\mathcal{M}=\max\left\{1,
\frac{2(\hat{c}+2C_{\rm dec}C_\epsilon)}{c}\right\}.
\]
The approximation in Theorem~\ref{thm2.6} extends this
estimate to energy solutions, since the approximate energies
converge uniformly on every finite time interval.

Letting $T\to\infty$ in \eqref{eq53} yields
\[
\int_s^\infty E(t)\dd t\le \mathcal{M}E(s),
\quad s\ge0.
\]
Applying Lemma~\ref{lem4.2} to $E(s+\cdot)$ yields
\eqref{eq54} for $t-s\ge \mathcal{M}$.
For $0\le t-s<\mathcal{M}$, monotonicity shows
\[
E(t)\le E(s)
\le \exp\!\left(1-\frac{t-s}{\mathcal{M}}\right)E(s).
\]
Hence \eqref{eq54} holds for all $t\ge s\ge0$.
\end{proof}
\subsection{Zero-flux boundary condition}
\label{sec:zero-flux}
We now consider the zero-flux case, whose weight assumptions
include classical strong degeneracy and the critical weight
$\widetilde a(x)=x^2$.
For $P(x,t)=\int_0^x\theta(y,t)\dd y$, the zero-flux condition
gives $P_t=a\theta_x-\kappa v$.
The remaining boundary term
$u_x(1,t)\int_0^1\theta\dd x$ is controlled by a mechanical
multiplier, yielding a strict Lyapunov inequality without
time averaging.

\begin{proposition}[A strict Lyapunov functional]
\label{prop4.6}
Assume Hypotheses~\ref{hyp2.2}-\ref{hybound}, in the zero-flux
case, and let $\kappa>0$.
There exist an instantaneous quadratic functional $\mathcal L_2$
and constants $c_1,c_2,c>0$, depending only on
$\kappa,M_1,M_2,\widetilde{M}$, such that every strong solution satisfies
\begin{equation}\label{eq55}
c_1E(t)\le\mathcal L_2(t)\le c_2E(t),
\quad
\mathcal L_2'(t)\le-cE(t).
\end{equation}
Consequently, every energy solution satisfies
\begin{equation}\label{eq56}
E(t)\le\frac{c_2}{c_1}
\exp\!\left(-\frac{c}{c_2}(t-s)\right)E(s),
\quad t\ge s\ge0.
\end{equation}
\end{proposition}

\begin{proof}
\emph{Step 1: multiplier identities.}
Let
\[
P(x,t)=\int_0^x\theta(y,t)\dd y,
\quad
\overline\theta(t)=\int_0^1\theta(y,t)\dd y.
\]
Then $P_x=\theta$, $P(0,t)=0$,
$P(1,t)=\overline\theta$ and $\|P\|\le\|\theta\|$ by Cauchy–Schwarz.
Integrating the thermal equation, using
$(a\theta_x)(0,t)=0$ and $v(0,t)=0$ gives
$
P_t=a\theta_x-\kappa v.
$

Define
\[
F=\int_0^1vP\dd x,\quad
G=\int_0^1uv\dd x,\quad
H=\int_0^1xu_xv\dd x.
\]
By Lemma~\ref{lem4.1}, the product rule and integration
by parts apply. Using $v_t=(u_x-\kappa\theta)_x$,
$P_t=a\theta_x-\kappa v$ and $P_x=\theta$, integration
by parts yields
\begin{equation}\label{eq57-F}
\begin{split}
F'
&=\int_0^1\bigl((u_x-\kappa\theta)_xP
  +v(a\theta_x-\kappa v)\bigr)\dd x\\
&=u_x(1,t)\overline\theta-\int_0^1u_x\theta\dd x
  +\kappa\|\theta\|^2+\int_0^1a\theta_xv\dd x
  -\kappa\|v\|^2,
\end{split}
\end{equation}
where $P(0,t)=0$, $P(1,t)=\overline\theta$ and
$\theta(1,t)=0$.

Similarly, using $u_t=v$ and the zero boundary values of $u$,
\begin{equation}\label{eq57-G}
\begin{aligned}
G'
=\|v\|^2+\int_0^1u(u_x-\kappa\theta)_x\dd x 
=\|v\|^2-\|u_x\|^2+\kappa\int_0^1u_x\theta\dd x.
\end{aligned}
\end{equation}

For $H$, local integration by parts on $(x,1)$ yields
\[
\begin{aligned}
\frac{\dd}{\dd t}\int_x^1yu_xv\dd y
&=\int_x^1y\bigl(v_xv+u_xu_{xx}
  -\kappa u_x\theta_x\bigr)\dd y\\
&=\frac12[y(v^2+u_x^2)]_{y=x}^{y=1}
  -\frac12\int_x^1(v^2+u_x^2)\dd y
  -\kappa\int_x^1y\theta_xu_x\dd y.
\end{aligned}
\]
The lower boundary terms vanish uniformly by \eqref{eq35},
and \eqref{eq20} gives
\[
\int_0^1|y\theta_xu_x|\dd y
\le\sqrt{M_1D}\,\|u_x\|.
\]
Passing to the limit $x\downarrow0$ in the time-integrated
identity and using $v(1,t)=0$, we obtain
\begin{equation}\label{eq57-H}
H'=\frac12u_x(1,t)^2
-\frac12(\|v\|^2+\|u_x\|^2)
-\kappa\int_0^1x\theta_xu_x\dd x.
\end{equation}

\emph{Step 2: the Lyapunov functional.}
Define
\[
\mathcal L_2=NE+F+\frac{\kappa}{4}G-\frac{\kappa}{4}H.
\]
Since $|F|,|G|,|H|\le E$, choosing
$N\ge2(1+\kappa/2)$ gives
\begin{equation}\label{eq58}
\frac N2E\le\mathcal L_2\le\frac{3N}{2}E.
\end{equation}
Using $E'=-D$ and \eqref{eq57-F}-\eqref{eq57-H}, we obtain
\begin{equation*}
\begin{aligned}
\mathcal L_2'&=-ND-\frac{5\kappa}{8}\|v\|^2
-\frac{\kappa}{8}\|u_x\|^2-\frac{\kappa}{8}u_x(1,t)^2
+u_x(1,t)\overline\theta\\
&\quad+\left(\frac{\kappa^2}{4}-1\right)
\int_0^1u_x\theta\dd x
+\kappa\|\theta\|^2+\int_0^1a\theta_xv\dd x
+\frac{\kappa^2}{4}\int_0^1x\theta_xu_x\dd x.
\end{aligned}
\end{equation*}
By Young's inequality and \eqref{eq20}, one has 
\[
\begin{aligned}
&|u_x(1,t)\overline\theta|
+\left|\left(\frac{\kappa^2}{4}-1\right)
\int_0^1u_x\theta\dd x\right|
+\left|\int_0^1a\theta_xv\dd x\right|
+\frac{\kappa^2}{4}\left|\int_0^1x\theta_xu_x\dd x\right|\\
&\quad\le
\frac{\kappa}{16}u_x(1,t)^2
+\frac{4M_2}{\kappa}D
+\frac{\kappa}{32}\|u_x\|^2
+\frac{8M_2}{\kappa}
\left(\frac{\kappa^2}{4}-1\right)^2D\\
&\qquad 
+\frac{\kappa}{16}\|v\|^2+\frac{4\widetilde{M}}{\kappa}D
+\frac{\kappa}{32}\|u_x\|^2+\frac{\kappa^3M_1}{2}D.
\end{aligned}
\]
Thus, setting
\begin{equation}\label{eq59}
C_D=\frac{4M_2}{\kappa}
+\frac{8M_2}{\kappa}\left(\frac{\kappa^2}{4}-1\right)^2
+\frac{4\widetilde{M}}{\kappa}+\frac{\kappa^3M_1}{2}
+\kappa M_2,
\end{equation}
we find
\[
\mathcal L_2'
\le-\frac{9\kappa}{16}\|v\|^2
-\frac{\kappa}{16}\|u_x\|^2
-\frac{\kappa}{16}u_x(1,t)^2-(N-C_D)D.
\]
Choose
\begin{equation}\label{eq60}
N\ge\max\{2(1+\kappa/2),C_D+1\},
\quad
c=\min\{\kappa/8,2/M_2\}.
\end{equation}
Since $\|\theta\|^2\le M_2 D$, it follows that
\[
\mathcal L_2'
\le-\frac{\kappa}{16}
\bigl(\|v\|^2+\|u_x\|^2\bigr)
-\frac1{M_2}\|\theta\|^2
\le-cE.
\]
Together with \eqref{eq58}, this proves
\eqref{eq55} with $c_1=N/2$ and $c_2=3N/2$.

\emph{Step 3: decay and energy data.}
Since $c_1E\le\mathcal L_2\le c_2E$ and
$\mathcal L_2'\le-cE$, we have
\[
\mathcal L_2'+\frac{c}{c_2}\mathcal L_2\le0.
\]
Integration and the equivalence with $E$ yield \eqref{eq56}.

The approximation in Theorem~\ref{thm2.6} gives uniform
convergence of the energies on compact time intervals, so this
estimate extends to energy solutions with the same constants.

\end{proof}

\begin{proof}[Proof of Theorem~\ref{thm2.7}]
The conclusion follows from
Propositions~\ref{prop4.5} and \ref{prop4.6}.
\end{proof}

\section{The sharp threshold for power conductivities}
\label{sec:powers}
Throughout this section, let $a(x)=x^\alpha$, $\alpha>0$, and
fix $\kappa>0$. We impose the thermal Dirichlet condition at zero
for $0<\alpha<1$ and the natural zero-flux condition for $\alpha\ge1$.
Define
\begin{equation*}
\begin{gathered}
V_\alpha=
\begin{cases}
\{z\in L^2(0,1)\cap AC([0,1]):
x^{\alpha/2}z_x\in L^2(0,1),\ z(0)=z(1)=0\},
&0<\alpha<1,\\
\{z\in L^2(0,1)\cap AC_{\rm loc}((0,1]):
x^{\alpha/2}z_x\in L^2(0,1),\ z(1)=0\},
&\alpha\ge1,
\end{cases}
\end{gathered}
\end{equation*}
equipped with the norm 
\begin{equation}
\label{eq61}
\|z\|_{V_\alpha}^2=\|z\|^2+\|x^{\alpha/2}z_x\|^2.
\end{equation}
For $0<\alpha<1$, the integrability of $\frac1{x^\alpha}$
ensures a well-defined trace at zero.
For $\alpha\ge1$, no temperature trace at zero is imposed
in $V_\alpha$.

Set
\begin{equation}
\label{119}
\mathcal A_\alpha(u,v,\theta)
=\bigl(v,(u_x-\kappa\theta)_x,
(x^\alpha\theta_x)_x-\kappa v_x\bigr),
\end{equation}
with 
\begin{equation}\label{eq62}
\begin{aligned}
\mathcal D_\alpha=\bigl\{(u,v,\theta)\in\mathcal H: 
v\in H_0^1(0,1), \theta\in V_\alpha, u_x-\kappa\theta,\ x^\alpha\theta_x\in H^1(0,1),(x^\alpha\theta_x)(0)=0
\ \text{if }\alpha\ge1\bigr\},
\end{aligned}
\end{equation}
These conditions ensure $\mathcal A_\alpha U\in\mathcal H$.
For $\alpha\ge1$, the zero-flux trace follows automatically from
$x^\alpha\theta_x\in H^1(0,1)$ and
\[
\int_0^1\frac{|x^\alpha\theta_x|^2}{x^\alpha}\dd x<\infty,
\quad
\int_0^1\frac{1}{x^{\alpha}}\dd x=\infty,
\]
by the argument of Lemma~\ref{lem3.2}. We retain it explicitly
to identify the boundary realization.

The boundary threshold $\alpha=1$ differs from the decay threshold
$\alpha=2$, since
\[
\sup_{0<x\le1}\frac{x^2}{x^\alpha}
=
\begin{cases}
1,&0<\alpha\le2,\\
+\infty,&\alpha>2.
\end{cases}
\]
The well-posedness result applies for every $\alpha>0$,
since Hypotheses~\ref{hyp2.2} and \ref{hyp2.3} hold with
$\widetilde a=a$ and $m=M=1$.
Hypothesis \ref{hybound}, and hence the preceding
decay theorem, applies for $0<\alpha\le2$.
For $\alpha>2$, the full norm in \eqref{eq61} still gives
well-posedness, while a concentration argument will exclude
uniform decay.

\begin{lemma}[Well-posedness without a weighted Poincare estimate]
\label{lem5.1}
Let $\alpha>2$ and $\kappa>0$. The operator
$(\mathcal A_\alpha,\mathcal D_\alpha)$ is densely defined and
maximal dissipative on $\mathcal H$, and generates a contraction
semigroup $S_\alpha(t)$.
For every $U_s\in\mathcal H$, the energy solution
$U(t)=S_\alpha(t-s)U_s$ satisfies
\begin{equation}\label{eq63}
\begin{gathered}
U\in C([s,T];\mathcal H),\quad
\theta\in L^2(s,T;V_\alpha),\\
E(t)+\int_s^t\|x^{\alpha/2}\theta_x(r)\|^2\dd r=E(s),
\quad s\le t\le T<\infty.
\end{gathered}
\end{equation}
If $U_s\in\mathcal D_\alpha$, then
\[
U\in C([s,T];\mathcal D_\alpha)\cap C^1([s,T];\mathcal H),
\]
where $\mathcal D_\alpha$ carries its graph norm.
\end{lemma}
\begin{proof}
For $a=\widetilde a=x^\alpha$, Hypothesis~\ref{hyp2.2} holds
with $m=M=1$ and the zero-flux realization, and
Hypothesis~\ref{hyp2.3} is automatic.
Moreover, $V=V_\alpha$, $\mathcal D(t)=\mathcal D_\alpha$
and $\mathcal A(t)=\mathcal A_\alpha$.
Proposition~\ref{prop3.5} gives maximal dissipativity and the
contraction semigroup. By uniqueness, its trajectories coincide
with the solutions of Theorem~\ref{thm2.6}, which gives
\eqref{eq63}, the regularity for domain data and the natural
zero-flux condition for energy data.
\end{proof}

The preceding lemma supplies the supercritical evolution independently
of the weighted lower bound. We can now combine the positive decay
result with a concentration argument to prove both parts of the threshold.

\begin{proof}[Proof of Theorem~\ref{thm2.9}]
We consider two cases.

\emph{Case 1: The range $0<\alpha\le2$.}
For $0<\alpha<1$,
\[
\int_0^1\frac{1}{x^{\alpha}}\dd x=\frac1{1-\alpha}<\infty,
\quad
\sup_{0<x\le1}x^{2-\alpha}=1,
\]
so the Dirichlet assumptions in Hypotheses~\ref{hyp2.2}
and  \ref{hybound} hold. 
For $1\le\alpha\le2$, the reciprocal integral diverges and 
the same supremum equals $1$, so the zero-flux assumptions
in Hypotheses~\ref{hyp2.2} and  \ref{hybound}
hold, including the critical case $\alpha=2$. 
Since $a$ is autonomous,
\[
\widetilde a=a,\quad m=M=1,\quad
M_2=4,\quad \widetilde{M}=M_1=1,
\]
and Hypothesis~\ref{hyp2.3} is automatic.
Proposition~\ref{prop3.5} and Theorem~\ref{thm2.7}
therefore give a contraction semigroup and exponential stability,
with decay constants independent of $\alpha$ in each boundary regime.
Taking the maximum of the prefactors and the minimum of the decay
rates gives a uniform exponential estimate for $0<\alpha\le2$.

\emph{Case 2: The range $\alpha>2$.}
By Lemma~\ref{lem5.1}, $\mathcal A_\alpha$ generates a contraction
semigroup $S_\alpha(t)$ on $\mathcal H$.
Choose a real function $\varphi\in C_c^\infty(1,2)$ with
$\|\varphi\|_{L^2(1,2)}=1$.
For $0<\varepsilon<1/2$, define
\[
\begin{aligned}
\theta_\varepsilon(x)
=\frac1{\sqrt\varepsilon}\varphi(\frac{x}{\varepsilon}),
\quad v_\varepsilon=0,\quad
u_\varepsilon(x)
=\kappa\left(
\int_0^x\theta_\varepsilon(y)\dd y
-x\int_0^1\theta_\varepsilon(y)\dd y
\right),
\end{aligned}
\]
and set
\begin{equation}\label{eq64}
U_\varepsilon=(u_\varepsilon,0,\theta_\varepsilon).
\end{equation}
The normalization preserves the $L^2$ norm of the temperature,
while its support concentrates near the degenerate endpoint.
Indeed, changing variables $x=\varepsilon y$ produces
\[
\|\theta_\varepsilon\|^2=1,
\quad
\int_0^1\theta_\varepsilon(x)\dd x
=\sqrt\varepsilon\int_1^2\varphi(y)\dd y.
\]
The displacement satisfies
\[
u_\varepsilon(0)=u_\varepsilon(1)=0,
\quad
(u_\varepsilon)_x-\kappa\theta_\varepsilon
=-\kappa\int_0^1\theta_\varepsilon(y)\dd y.
\]
Thus its stress is spatially constant.
Since the temperature and flux are smooth and supported in
$(\varepsilon,2\varepsilon)$, we have
$U_\varepsilon\in\mathcal D_\alpha$ and
$
\mathcal A_\alpha U_\varepsilon
=\left(0,0,(x^\alpha(\theta_\varepsilon)_x)_x\right).
$

The constant-stress construction yields
\[
\|(u_\varepsilon)_x\|^2
=\kappa^2\left(
\|\theta_\varepsilon\|^2
-\left|\int_0^1\theta_\varepsilon\dd x\right|^2
\right).
\]
Moreover,
$
(x^\alpha(\theta_\varepsilon)_x)_x(\varepsilon y)
=\varepsilon^{\alpha-5/2}(y^\alpha\varphi'(y))'.
$ 
Thus, together with $\alpha>2$, as $\varepsilon\downarrow0$,
\begin{equation}\label{eq65}
\begin{aligned}
\|U_\varepsilon\|_{\mathcal H}^2
&=1+\kappa^2-\kappa^2\varepsilon
\left|\int_1^2\varphi(y)\dd y\right|^2
\to1+\kappa^2,\\
\|\mathcal A_\alpha U_\varepsilon\|_{\mathcal H}^2
&=\varepsilon^{2\alpha-4}
\int_1^2|(y^\alpha\varphi'(y))'|^2\dd y
\to0.
\end{aligned}
\end{equation}
Hence these states have nonvanishing energy and a vanishing
evolution residual.

For $U_\varepsilon\in\mathcal D_\alpha$, the semigroup identity
$
S_\alpha(t)U_\varepsilon-U_\varepsilon
=\int_0^tS_\alpha(r)\mathcal A_\alpha U_\varepsilon\dd r
$
and contraction imply
\begin{equation}\label{eq66}
\|S_\alpha(t)U_\varepsilon-U_\varepsilon\|_{\mathcal H}
\le\int_0^t
\|S_\alpha(r)\mathcal A_\alpha U_\varepsilon\|_{\mathcal H}\dd r
\le\int_0^t
\|\mathcal A_\alpha U_\varepsilon\|_{\mathcal H}\dd r
=t\|\mathcal A_\alpha U_\varepsilon\|_{\mathcal H}.
\end{equation}
The same estimate also quantifies the persistence of these states.
Set
\[
C_{\alpha,\varphi}
=\left(\int_1^2|(y^\alpha\varphi'(y))'|^2\dd y\right)^{1/2}>0.
\]
Since $\|U_\varepsilon\|_{\mathcal H}\ge
\|\theta_\varepsilon\|=1$, equations \eqref{eq65}--\eqref{eq66}
give, for every $t\ge0$,
\[
\frac{\|S_\alpha(t)U_\varepsilon-U_\varepsilon\|_{\mathcal H}}
{\|U_\varepsilon\|_{\mathcal H}}
\le C_{\alpha,\varphi}t\varepsilon^{\alpha-2}.
\]
Writing $E_\varepsilon(t)=\frac12
\|S_\alpha(t)U_\varepsilon\|_{\mathcal H}^2$, the reverse
triangle inequality therefore yields
\[
E_\varepsilon(t)\ge\frac14E_\varepsilon(0),
\quad
0\le t\le
\frac{1}{2C_{\alpha,\varphi}\varepsilon^{\alpha-2}}.
\]
For fixed $\alpha>2$, this time interval becomes arbitrarily
long as $\varepsilon\downarrow0$.

Thus, for every fixed $t\ge0$, inequality $
\|S_\alpha(t)U_\varepsilon\|_{\mathcal H}
\ge
\|U_\varepsilon\|_{\mathcal H}
-\|S_\alpha(t)U_\varepsilon-U_\varepsilon\|_{\mathcal H}$ yields
\[
\begin{aligned}
1
\ge\|S_\alpha(t)\|_{\mathcal L(\mathcal H)}=
\sup_{W\ne0}
\frac{\|S_\alpha(t)W\|_{\mathcal H}}{\|W\|_{\mathcal H}}
\ge
\frac{\|S_\alpha(t)U_\varepsilon\|_{\mathcal H}}
{\|U_\varepsilon\|_{\mathcal H}} \ge
1-t\frac{\|\mathcal A_\alpha U_\varepsilon\|_{\mathcal H}}
{\|U_\varepsilon\|_{\mathcal H}}
\to1\quad\text{as }\varepsilon\downarrow0.
\end{aligned}
\]
Therefore
$\|S_\alpha(t)\|_{\mathcal L(\mathcal H)}=1$ for every $t\ge0$,
proving \eqref{eq17}.

If $E(t)\le r(t)E(0)$ held for all energy initial data, since $
U(t)=S_\alpha(t)U_0,$ $
E(t)=\frac12\|S_\alpha(t)U_0\|_{\mathcal H}^2,$ we have 
\[
1=\|S_\alpha(t)\|_{\mathcal L(\mathcal H)}^2
=\sup_{U_0\ne0}
\frac{\|S_\alpha(t)U_0\|_{\mathcal H}^2}
{\|U_0\|_{\mathcal H}^2}
\le r(t).
\]
Hence no such uniform decay factor can tend to zero as
$t\to\infty$.
\end{proof}

\begin{corollary}[No uniform decay from the graph domain]
\label{cor5.2}
Let $\alpha>2$ and equip $\mathcal D_\alpha$ with
\[
\|U\|_{\mathcal D_\alpha}
=\|U\|_{\mathcal H}+\|\mathcal A_\alpha U\|_{\mathcal H}.
\]
Then
\begin{equation}\label{eq67}
\|S_\alpha(t)\|_{\mathcal L(\mathcal D_\alpha,\mathcal H)}=1,
\quad t\ge0.
\end{equation}
Consequently, no function $r(t)\to0$ as $t\to\infty$ satisfies
\[
\|S_\alpha(t)U\|_{\mathcal H}^2
\le r(t)\|U\|_{\mathcal D_\alpha}^2,
\quad U\in\mathcal D_\alpha.
\]
The same nonexistence holds for any equivalent graph norm.
\end{corollary}

\begin{proof}
By contraction and the definition of the graph norm,
$
\|S_\alpha(t)U\|_{\mathcal H}
\le\|U\|_{\mathcal H}
\le\|U\|_{\mathcal D_\alpha},
$
so
$\|S_\alpha(t)\|_{\mathcal L(\mathcal D_\alpha,\mathcal H)}
\le1$.
For the states in \eqref{eq64}, the reverse triangle inequality
and \eqref{eq66} yield
\[
\begin{aligned}
\|S_\alpha(t)\|_{\mathcal L(\mathcal D_\alpha,\mathcal H)}
\ge
\frac{\|S_\alpha(t)U_\varepsilon\|_{\mathcal H}}
{\|U_\varepsilon\|_{\mathcal D_\alpha}}
\ge
\frac{\|U_\varepsilon\|_{\mathcal H}
-t\|\mathcal A_\alpha U_\varepsilon\|_{\mathcal H}}
{\|U_\varepsilon\|_{\mathcal H}
+\|\mathcal A_\alpha U_\varepsilon\|_{\mathcal H}}
\to1\quad\text{as }\varepsilon\downarrow0
\end{aligned}
\]
for every fixed $t\ge0$, since \eqref{eq65} gives
$\|U_\varepsilon\|_{\mathcal H}\to\sqrt{1+\kappa^2}$
and $\|\mathcal A_\alpha U_\varepsilon\|_{\mathcal H}\to0$.
Together with the upper bound, this proves \eqref{eq67}.

If the stated decay estimate held for all
$U\in\mathcal D_\alpha$, taking the supremum over nonzero
initial data would give
\[
1=\|S_\alpha(t)\|_{\mathcal L(\mathcal D_\alpha,\mathcal H)}^2
\le r(t),
\]
which excludes $r(t)\to0$ as $t\to\infty$.
For any equivalent graph norm $\|\cdot\|_*$, there exists
$C>0$ such that
$\|U\|_*\le C\|U\|_{\mathcal D_\alpha}$.
A decay estimate with $\|\cdot\|_*$ would therefore imply
the same estimate with the original graph norm and factor
$C^2r(t)\to0$ as $t\to\infty$, giving the same contradiction.
\end{proof}

\begin{proposition}[The quadratic coefficient]
\label{prop5.3}
For $a(x)=x^2$ with the natural zero-flux realization,
the embedding $V_2\hookrightarrow L^2(0,1)$ is not compact,
but
\begin{equation}\label{eq68}
\|z\|^2\le4\|xz_x\|^2,\quad z\in V_2,
\end{equation}
and the coupled system is exponentially stable.
\end{proposition}

\begin{proof}
Estimate \eqref{eq68} follows from Lemma~\ref{lem3.1}
with $a=x^2$.
For the temperature profiles in \eqref{eq64}, the change of
variables $x=\varepsilon y$ shows
$
\|\theta_\varepsilon\|^2=1,$ $
\|x\theta_{\varepsilon,x}\|^2
=\int_1^2y^2|\varphi'(y)|^2\dd y.
$
Hence
$
\|\theta_\varepsilon\|_{V_2}^2
=1+\int_1^2y^2|\varphi'(y)|^2\dd y,
$
which is independent of $\varepsilon$.
Taking $\varepsilon_n=\frac1{3^n}$, we have
$
\operatorname{supp}\theta_{\varepsilon_n}
\subset
\left(\frac1{3^n},\frac2{3^n}\right).
$
These intervals are pairwise disjoint, so
\[
\|\theta_{\varepsilon_n}-\theta_{\varepsilon_m}\|^2
=\|\theta_{\varepsilon_n}\|^2
+\|\theta_{\varepsilon_m}\|^2
=2,\quad n\ne m.
\]
Thus this bounded sequence in $V_2$ has no convergent
subsequence in $L^2(0,1)$, proving that the embedding is
not compact.

For $a=\widetilde a=x^2$, the constants satisfy
$M_2=4$ and $\widetilde M=M_1=1$, so exponential stability
follows from Theorem~\ref{thm2.9} with $\alpha=2$.
Thus the weighted coercivity estimate provides decay even
though the thermal embedding is not compact.
\end{proof}

\begin{remark}[An explicit critical decay bound]
\label{rem5.5}
For $a=x^2$ and $\kappa=1$, we have
$M_2=4$ and $\widetilde{M}=M_1=1$, hence
$C_D=85/2$ by \eqref{eq59}.
Taking $N=44$ and $c=1/8$ in
\eqref{eq60} gives
$c_1=22$ and $c_2=66$. Thus
\eqref{eq56} yields
\[
E(t)\le3e^{-(t-s)/528}E(s),\quad t\ge s\ge0.
\]
This bound is not claimed to be optimal.
\end{remark}

\section{Examples and scope of the assumptions}
\label{sec:extensions}

The first three examples separate the spatial bounds that determine
the decay constants from the ratio regularity used to construct the
evolution. The last example shows that the Dirichlet theorem extends
beyond the classical degeneracy index. Each example admits global
energy solutions satisfying \eqref{eq16} for every fixed
$\kappa>0$.

\begin{example}[A critical conductivity with rapid temporal oscillations]
\label{ex6.1}
Let $0<\varepsilon<1$ and impose the natural zero-flux condition.
The conductivity
\begin{equation*}
a(x,t)=x^2\bigl(1+\varepsilon\sin(t^2)\sin(\log x)\bigr)
\quad\text{for }x>0,\quad a(0,t)=0,
\end{equation*}
satisfies Hypotheses~\ref{hyp2.2}-\ref{hybound}. 
Its classical index exceeds $2$, its ratio has no endpoint limit
when $\sin(t^2)\ne0$, and $|a_t/a|$ is unbounded on the time
half-line. The decay constants depend only on $\varepsilon$ and
$\kappa$.
\end{example}

\begin{proof}[Verification]
Here $\widetilde a=x^2$ and $1-\varepsilon\le a/\widetilde a=1+\varepsilon\sin(t^2)\sin(\log x)\le1+\varepsilon$. Direct differentiation gives
\begin{equation}\label{eq70}
\begin{aligned}
&\left(\frac {a}{\widetilde a}\right)_t=2\varepsilon t\cos(t^2)\sin(\log x),\quad
x\left(\frac {a}{\widetilde a}\right)_x=\varepsilon\sin(t^2)\cos(\log x),\\
&x\left(\frac {a}{\widetilde a}\right)_{tx}=2\varepsilon t\cos(t^2)\cos(\log x),\quad
a_x=2x\frac {a}{\widetilde a}+\varepsilon x\sin(t^2)\cos(\log x).
\end{aligned}
\end{equation}
These formulas verify \eqref{eq8} on every finite interval
and give $a(\cdot,t)\in C^1[0,1]$ with $a_x(0,t)=0$.
The weight conditions hold with
\[
M_2=\frac4{1-\varepsilon},\quad
M_1=\frac1{1-\varepsilon},\quad \widetilde{M}=1+\varepsilon.
\]
At $t=\sqrt{\pi/2}$ and $x=\frac{1}{e^{2\pi n}}$, \eqref{eq70}
gives $x|a_x|/a=2+\varepsilon$. For $\sin(t^2)\ne0$, the
oscillation of $\sin(\log x)$ prevents a limit of $\frac {a}{\widetilde a}$ at zero.
Finally, at $t_n=\sqrt{2\pi n}$ and $x_*=\frac{1}{e^{\pi/2}}$,
$|a_t/a|(x_*,t_n)=2\varepsilon t_n\to\infty$.
By continuity away from zero, the spacetime essential supremum is
also infinite, whereas \eqref{eq10} holds locally.
Theorems~\ref{thm2.6} and \ref{thm2.7} give the conclusion.
\end{proof}

\begin{example}[A weakly degenerate ratio with unbounded derivative]
\label{ex6.2}
Let $0<\alpha<1$, $1/2<\beta<1$, and $0<\varepsilon<1$ satisfy
\begin{equation}\label{eq71}
\alpha+\frac{\varepsilon\beta}{1-\varepsilon}<1.
\end{equation}
Under the Dirichlet condition, set
\begin{equation*}
a(x,t)=x^\alpha(1+\varepsilon\sin t\,x^\beta).
\end{equation*}
This coefficient is classically \textup{(WD)}. Its ratio satisfies
{\eqref{eq7}} {in the Dirichlet case}, although
${\partial_x(a(\cdot,t)/\widetilde a)\notin L^\infty(0,1)}$ when $\sin t\ne0$.
\end{example}

\begin{proof}[Verification]
Here $\widetilde a=x^\alpha$ and 
$1-\varepsilon\le a/\widetilde a=1+\varepsilon\sin t\,x^\beta\le1+\varepsilon$. Moreover,
\[
\frac{x|a_x|}{a}\le\alpha+\frac{\varepsilon\beta}{1-\varepsilon}<1.
\]
Thus \eqref{eq71} gives classical weak
degeneracy. Since $x^{\beta-1}\in L^2(0,1)$ for $\beta>1/2$,
$\frac {a}{\widetilde a}$ and $\left(\frac {a}{\widetilde a}\right)_t=\varepsilon\cos t\,x^\beta$ satisfy
\eqref{eq7}. But
$\left(\frac {a}{\widetilde a}\right)_x=\varepsilon\beta\sin t\,x^{\beta-1}$ is unbounded when
$\sin t\ne0$, because $\beta<1$.
Corollary~\ref{cor2.8} applies.
\end{proof}

\begin{example}[A strongly degenerate ratio without an endpoint limit]
\label{ex6.3}
Let $1<\alpha<2$ and $0<\varepsilon<1$ satisfy
\begin{equation*}
\alpha+\frac{\varepsilon}{1-\varepsilon}<2.
\end{equation*}
Under the natural zero-flux condition, set
\begin{equation*}
a(x,t)=x^\alpha\bigl(1+\varepsilon\sin t\sin(\log x)\bigr)
\quad\text{for }x>0,\quad a(0,t)=0.
\end{equation*}
This coefficient is classically \textup{(SD)} and satisfies
{\eqref{eq8}} {in the zero-flux case}. When $\sin t\ne0$, its ratio
has no continuous extension to zero and does not belong to $H^1(0,1)$.
\end{example}

\begin{proof}[Verification]
Here $\widetilde a=x^\alpha$ and 
$1-\varepsilon\le a/\widetilde a=1+\varepsilon\sin t\sin(\log x)\le 1+\varepsilon$. The derivative formulas show
\[
\frac{x|a_x|}{a}\le\alpha+\frac{\varepsilon}{1-\varepsilon}<2,
\quad a_x=O(x^{\alpha-1})\quad\text{for }x\downarrow0.
\]
Thus $a(\cdot,t)\in C^1[0,1]$ with $a_x(0,t)=0$, and
$1<\alpha\le K<2$ since the index at $t=0$ equals $\alpha$.
The verification of \eqref{eq8} is the same as in
Example~\ref{ex6.1}, with $\sin(t^2)$ replaced
by $\sin t$. Corollary~\ref{cor2.8} therefore applies.
When $\sin t\ne0$, oscillation prevents an endpoint limit, and
\[
\int_0^1\left|\left(\frac {a}{\widetilde a}\right)_x\right|^2\dd x
=\varepsilon^2\sin^2t\int_0^\infty e^y\cos^2y\dd y=\infty.
\]
Hence $\frac {a}{\widetilde a}\notin H^1(0,1)$ at those times.
\end{proof}

\begin{example}[{An oscillating coefficient with infinite index}]
\label{ex6.4}
For $0<\alpha<1$, consider the autonomous coefficient
\begin{equation*}
a(x)=x^\alpha\bigl(2+\sin(1/x)\bigr)\quad\text{for }x>0,
\quad a(0)=0,
\end{equation*}
with the Dirichlet thermal condition. It satisfies the hypotheses
of Theorems~\ref{thm2.6} and \ref{thm2.7}, although
its classical index is infinite.
\end{example}

\begin{proof}[Verification]
The coefficient belongs to $C[0,1]\cap C^\infty((0,1])$ and satisfies
$x^\alpha\le a(x)\le 3x^\alpha$. Hence 
\[
\int_0^1\frac{dx}{a(x)}
\le\int_0^1x^{-\alpha}\,dx
=\frac1{1-\alpha}<\infty,\quad 
\sup_{0<x\le1}\frac{x^2}{a(x)}
\le \sup_{0<x\le1}x^{2-\alpha}\le1.
\]
With $\widetilde a=a$, the ratio is identically one, so all coefficient
hypotheses hold. At $x_n=\frac{1}{2\pi n}$,
\[
\frac{x_na'(x_n)}{a(x_n)}=\alpha-\pi n\to-\infty\quad\text{as }n\to \infty,
\]
which proves that the classical index is infinite.
This coefficient also satisfies the weakly degenerate assumptions and
$W\in L^\infty(0,1)$ in \cite{AHR}.
Its autonomous exponential stability for the sufficiently small values
of $\kappa$ considered there is already covered by their Theorem~3.3.
Here the example illustrates the distinction between the weight
assumptions and a finite classical index.
\end{proof}

\begin{remark}[Scope of the results]\label{rem6.5}
Hypotheses~\ref{hyp2.2} and \ref{hyp2.3} are sufficient
for the common-domain construction. The decay argument additionally
uses Hypothesis \ref{hybound}. No minimality is claimed. 
Uniform comparison with a fixed coefficient requires a
time-independent power exponent. Indeed, if $a(x,t)=x^{\alpha(t)}$,
then
$a(x,t)/a(x,0)=x^{\alpha(t)-\alpha(0)}$ tends to zero or infinity
as $x\downarrow0$ whenever $\alpha(t)\ne\alpha(0)$.
Discontinuous switching is also outside the time regularity in
{Hypothesis~\ref{hyp2.3}}.

The threshold $\alpha=2$ concerns uniform decay for the autonomous
power family with the stated boundary realizations. It does not give
a necessary and sufficient criterion for general conductivities, or
describe systems with time-dependent boundary realizations.
For $\alpha>2$,
\eqref{eq17} excludes uniform decay on the energy space,
and Corollary~\ref{cor5.2} excludes it even from the
generator graph domain into that space. These conclusions leave
strong stability of individual trajectories open. They do not address
decay on other spaces whose norms exclude the concentrating sequence
\eqref{eq64} from bounded sets.
\end{remark}

\section*{Declaration of competing interest}

The author declares that they have no known competing financial interests or personal relationships that could have appeared to influence the work reported in this paper.

\section*{Data availability}

Data sharing is not applicable to this article as no data sets were generated or analysed during the current study.

\section*{Declaration of generative AI and AI-assisted technologies in the manuscript preparation process}
During the preparation of this work, the author used GPT-6.1 Sol for textual revision and polishing. Discussions with this AI tool inspired some research ideas and assisted with the construction of illustrative examples. All core proof strategies, rigorous mathematical arguments and final validations were independently completed and checked by the author. The author reviewed and edited all relevant outputs and takes full responsibility for the content of the published article.

\end{document}